\documentclass{amsart}
\usepackage{amsmath, amssymb, xspace, epsfig, float, array}
\usepackage{latexsym}
\usepackage{alltt}
\usepackage{fancyhdr}
\usepackage{subfigure}
\usepackage{algorithm}
\usepackage{algorithmic}
\usepackage{graphicx}
\include{elsart}

\newtheorem{theorem}{Theorem}[section]
\newtheorem{corollary}{Corollary}[section]
\newtheorem{definition}{Definition}[section]
\newtheorem{remark}{Remark}[section]
\newtheorem{lemma}{Lemma}[section]
\newtheorem{proposition}{Proposition}[section]
\newtheorem*{theorem*}{Theorem}
\def\VV{{V}}

\def\CC{{C}}
\DeclareMathOperator{\gds}{{g_{\it d}sin}}
\DeclareMathOperator{\gdms}{{g_{{\it d}-1}sin}}
\DeclareMathOperator{\pds}{{p_{\it d}sin}}
\DeclareMathOperator{\pons}{{p_{1}sin}}
\DeclareMathOperator{\pdms}{{p_{{\it d}-1}sin}}

\def\WW{{W}}

\def\BB{{B}}
\def\SS{{S}}
\DeclareMathOperator{\dist}{dist}
\DeclareMathOperator{\sign}{sign}
\DeclareMathOperator{\diam}{diam}
\DeclareMathOperator{\Sp}{Sp}
\def\LL{{L}}
\DeclareMathOperator{\DIm}{dim}
\def\HH{{H}}
\def\RR{{\mathbb{R}}}
\newcommand{\reals}{\mathbb R}
\newcommand{\ints}{\mathbb Z}
\newcommand{\nats}{\mathbb N}

\def\Supp{{\mathrm{supp}(\mu)}}
\def\Pllp{{\mathrm{P}_{\mathrm{rll}}}}
\def\Cpoly{{\mathrm{C}_{\mathrm{poly}}}}
\def\Affn{{\mathrm{A}_{\mathrm{ffn}}}}
\def\Chull{{\mathrm{C}_{\mathrm{hull}}}}
\def\Cne{{\mathrm{C}_{\mathrm{one}}}}
\def\Cnetilde{{\mathrm{\tilde{C}_{one}}}}
\def\Ttilde{{\tilde{\Theta}}}
\def\Tbe{{\mathrm{T_{ube}}}}
\usepackage[top=1in,bottom=1in,left=1in,right=1in]{geometry}
\newcommand{\LV}[2]{\LL_{{#1}\setminus\{{#2}\}}}
\newcommand{\LVU}[3]{\LL_{{#1}\setminus{\{{#2},{#3}\}}}}
\newcommand{\LS}[1]{\LL_{#1}}
\newcommand{\TV}[3]{\theta\left({#3},\LV{#1}{#2}\right)}
\newcommand{\TVU}[4]{\theta\left({#4},\LVU{#1}{#2}{#3}\right)}

\newcommand{\AVU}[3]{\alpha\left(\LV{#1}{#2},\LV{#1}{#3}\right)}

\begin{document}

\title[$d$-dimensional $d$-Semimetrics and Simplex
Inequalities for High-Dimensional Sines]{On
$d$-dimensional $d$-Semimetrics and Simplex-Type Inequalities\\
for High-Dimensional Sine Functions}
\author{Gilad Lerman}
\address{Department of Mathematics, University of Minnesota,
        127 Vincent Hall, 206 Church St. S.E., Minneapolis, MN 55455
        USA.}
\email{lerman@math.umn.edu}

\author{J. Tyler Whitehouse}
\address{Department of Mathematics, University of Minnesota,
        127 Vincent Hall, 206 Church St. S.E., Minneapolis, MN 55455
        USA.}
\email{whit0933@umn.edu}
\thanks{The authors have been supported by NSF grant \#0612608}

\date{June 10, 2007}
%
%


\begin{abstract}
We show that high-dimensional analogues of the sine function (more
precisely, the $d$-dimensional polar sine
and the $d$-th root of the $d$-dimensional hypersine)
satisfy a simplex-type inequality
in a real pre-Hilbert space $\HH$. Adopting the language of Deza and
Rosenberg,
we say that these $d$-dimensional sine
functions are $d$-semimetrics.
We also establish geometric identities for both the $d$-dimensional
polar sine and the $d$-dimensional hypersine.
We then show that when $d=1$ the underlying functional
equation of the corresponding identity
characterizes a generalized sine function.
%
Finally, we
show that the $d$-dimensional
polar sine satisfies a relaxed simplex inequality of two controlling
terms ``with high probability''.

\end{abstract}
\maketitle


\section{Introduction}We establish some fundamental properties of high-dimensional sine
functions~\cite{eriksson}. In particular, we show that they
satisfy a simplex-type inequality, and thus according to the
terminology of Deza and Rosenberg~\cite{DR00} they are
$d$-semimetrics. We also demonstrate a related concentration
inequality. These properties are useful to some modern
investigations in harmonic analysis~\cite{LW-part1,LW-part2} and
applied mathematics~\cite{spectral_theory,spectral_applied}.

High-dimensional sine functions have been known for more than a
century. Euler~\cite{Euler} formulated the two-dimensional polar
sine (for tetrahedra) and D'Ovidio~\cite{Ovidio} generalized it to
higher dimensions. Joachimsthal~\cite{Joachimsthal} suggested the
two-dimensional hypersine (for tetrahedra) and
Barto\v{s}~\cite{bartos} extended it to simplices of any
dimension. Various authors have explored their properties and
applied them to a variety of problems (see e.g., \cite{eriksson},
\cite{V95}, \cite{LY98}, \cite{XL05}, \cite{polar_duality} and
references in there). For our purposes, we have slightly modified
the existing definitions, in particular we allow negative values
of these functions when the dimension of the ambient space is
$d+1$.

The two high-dimensional sine functions that we define here, $\pds$
and $\gds$, are exemplified in Figure~\ref{figure:def_sines} and are
described as follows.
For $v_1,\ldots,v_{d+1}, w \in \HH$ we take the parallelotope
through the points $v_1,\ldots,v_{d+1}, w$.
The function $|\pds_w(v_1,\ldots,v_{d+1})|$ is obtained by dividing
the ($d+1$)-volume of that parallelotope by the $d+1$ edge lengths
at the vertex $w$.
Similarly, we define  $|\gds_w(v_1,\ldots,v_{d+1})|$ to be the
($d+1$)-volume of the same parallelotope scaled by the $d$-th
roots of the $d$-volumes of its faces through the vertex $w$
(there are $d+1$ of these). That is,
$|\pds_w(v_1,\ldots,v_{d+1})|$ and $|\gds_w(v_1,\ldots,v_{d+1})|$
are the polar sine~\cite{eriksson} and the $d$-th root of the
hypersine~\cite{gsin_mathworld} of the simplex with vertices
$\{w,v_1,\ldots,v_{d+1}\}$ with respect to the vertex $w$. If
$\DIm(H) = d+1$, then we define $\pds_w(v_1,\ldots,v_{d+1})$ and
$\gds_w(v_1,\ldots,v_{d+1})$ by replacing the volume of the
parallelotope by the corresponding determinant
(precise definitions appear in
Subsection~\ref{subsection:sine-definition}). We often assume that
$w=0$ since the more general case can be obtained by a simple shift
(as expressed later in equations~\eqref{eq:gsin_invariant}
and~\eqref{eq:psin_invariant}).

We note that when $d=1$: $|\pons_0(v_1,v_2)| = |g_1 \sin_0(v_1,v_2)|
= |\sin(v_1,v_2)|$, where $|\sin(v_1,v_2)|$ denotes the absolute
value of the sine of the angle between $v_1 $ and  $v_2$.
Furthermore, regardless of the dimension of $\HH$, the following
triangle inequality holds
\begin{equation}\label{eq:regsine2} |\sin(v_1,v_2)| \leq |
\sin(v_1,u)| + |\sin(u,v_2)|,\textup{ for all } v_1, v_2 \in \HH
 \textup{ and } u \in \HH \setminus \{0\}.
\end{equation}The first part of this paper establishes high-dimensional analogues of
equation~(\ref{eq:regsine2}) for the functions $|\pds|$ and $|\gds|$
for $d>1$.

One motivation for our research is the interest in high-dimensional
versions of metrics
and $d$-way kernel methods in machine learning~\cite{Agarwal05,
Shashua06}. Deza and Rosenberg~\cite{DR00} have defined the notion
of a $d$-semimetric (or $n$-semimetric according to their
notation). If $d \in \nats$ and $E$ is a given set,  then the pair
$(E,f)$ is a $d$-semimetric if $f :E^{d+1} \mapsto [0,\infty)$ is
symmetric (invariant to permutations) and satisfies the following
simplex-type inequality:
\begin{equation}
\label{eq:simplex-type_deza} f(x_1,\ldots,x_{d+1}) \leq
\sum_{i=1}^{d+1} f(x_1,\ldots,x_{i-1},u,x_{i+1},\ldots,x_{d+1})\, \
\textup{ for all } \ x_1,\ldots,x_{d+1},u \in E.
\end{equation}
We also refer to $f$ itself as a $d$-semimetric with respect to $E$
or just a $d$-semimetric when the set E is clear.

The examples of $d$-semimetrics proposed by Deza and
Rosenberg~\cite{DR00} do not represent $d$-dimensional geometric
properties. They typically form $d$-semimetrics by averaging
non-negative functions that quantify lower order geometric
properties of a $d$-simplex (see~\cite[Fact~2]{DR00}). For example,
in order to form a $2$-semimetric on $\HH$ they average the pairwise
distances between three points to get the scaled perimeter of the
corresponding triangle, which is a one-dimensional quantity.

We provide here the following $d$-dimensional examples of
$d$-semimetrics.
\begin{theorem}\label{theorem:relaxedboth} If $\HH$ is a real pre-Hilbert
space, $d\in\mathbb{N}$, and $\DIm(\HH)\geq d+1$, then the
functions $|\pds_0|$ and $|\gds_0|$ are $d$-semimetrics with
respect to the set $\HH\setminus \{0\}$.
\end{theorem}
The above examples of $d$-semimetrics are $d$-dimensional in the
following sense: $|\pds_0(v_1,\ldots,v_{d+1})|$ and
$|\gds_0(v_1,\ldots,v_{d+1})|$ are zero if and only if the vectors
$v_1,\ldots,v_{d+1}$ are linearly dependent, and they are one (and
maximal) if and only if the vectors $v_1,\ldots,v_{d+1}$ are
mutually orthogonal.

%

Another motivation for our research is our interest in
high-dimensional generalizations of the Menger
curvature~\cite{M30,MMV96}. In a subsequent
work~\cite{LW-part1,LW-part2} we define a $d$-dimensional
Menger-type curvature for $d > 1$ via the polar sine, $|\pds|$,
and use it to characterize the smoothness of $d$-dimensional
Ahlfors regular measures (see Definition~\ref{def:ARmeasure}). Our
proof utilizes the fact that the polar sine satisfies ``a relaxed
simplex inequality of two controlling terms with high Ahlfors
probability''. We quantify this notion in a somewhat general
setting as follows.

For a symmetric function $f$ on $H^{d+1}$, an integer $p$, $1 \leq p
\leq d$, and a positive constant $C$, we say that $f$ satisfies a
{\em relaxed simplex inequality} of $p$ terms and constant $C$ if
\begin{equation}\label{eq:relaxed-intro}f(v_1,\ldots,v_{d+1})\leq
C\cdot\sum_{i=1}^p
f(v_1,\ldots,v_{i-1},u,v_{i+1},\ldots,v_{d+1}),\textup{ for all
}v_1,\ldots,v_{d+1}\in\HH\textup{ and
}u\in\HH\setminus\{0\}.\end{equation}
By the symmetry of $f$, any $p$ terms in the above sum will suffice
(e.g., replacing $\sum_{i=1}^p$ by $\sum_{i=d+2-p}^{d+1}$ in
equation~(\ref{eq:relaxed-intro})).

For any $\SS=\{v_1,\ldots,v_{d+1}\}\subseteq\HH$, $w\in\HH$, and
$C>0$, we let $U_C(\SS,w)$ be the set of vectors $u$ giving rise to
relaxed simplex inequalities of two terms and constant $C$ for
$|\pds_w|$, that is,
\begin{multline}\label{def:U-C}U_C(\SS,w)=\Big\{u\in\HH:\ |
\pds_w(v_1,\ldots,v_{d+1})|\leq
C\cdot\big(|\pds_w(v_1,\ldots,v_{i-1},u,v_{i+1},\ldots,v_{d+1})|\
+\\|\pds_w(v_1,\ldots,v_{j-1},u,v_{j+1},\ldots,v_{d+1})|\big),\textup{
for all }1\leq i<j\leq d+1\Big\}.\end{multline}
Using this notation, we claim that for any $d$-dimensional Ahlfors
regular measure $\mu$ on $\HH$, any sufficiently large constant $C$,
any set of vectors $\SS$ as above and any $w\in\Supp$, the event
$U_C(\SS,w)$ has high probability at any relevant ball in $\HH$,
where a probability at a ball is obtained by scaling the measure
$\mu$ by the measure of the ball.
We formulate this property more precisely and even more generally as
follows:
\begin{theorem}\label{theorem:relaxedar} If $\HH$ is a pre-Hilbert space,
$2 \leq d \in \mathbb{N}$,  $0< \epsilon < 1$, $\gamma\in\RR$ is
such that $d-1<\gamma\leq d$, $w\in\Supp$,
$\SS=\{v_1,\ldots,v_{d+1}\}\subseteq\HH$, and $\mu$ is a
$\gamma$-dimensional Ahlfors regular measure on $\HH$ with Ahlfors
regularity constant $\CC_{\mu}$,  then there exists a constant
$C_0\geq 1$ depending only on $C_{\mu}$, $\epsilon$, $\gamma$, and
$d$, such that for all $C\geq C_0$:
\begin{equation}
\label{eq:theorem2}
\frac{\mu\left(U_C(\SS,w)\cap\BB(w,r)\right)}{\mu\left(\BB(w,r)\right)}\geq
1-\epsilon, \ \textup{ for all } 0<r\leq\diam(\Supp).\end{equation}
\end{theorem}


The paper is organized as follows.  In
Section~\ref{section:notation} we present the main notation and
definitions as well as a few elementary properties of the
$d$-dimensional sine functions. In
Section~\ref{section:identities} we develop geometric identities
for $\pds_0$ and $\gds_0$ as well as characterize the solutions of
the corresponding functional equations when $d=1$.  In
Section~\ref{section:inequalities} we prove
Theorem~\ref{theorem:relaxedboth},
and in Section~\ref{section:relaxed} we prove
Theorem~\ref{theorem:relaxedar}.
Finally, we conclude
our research in Section~\ref{section:final} and discuss future
directions and open problems.

\section{Notation, Definitions, and Elementary
Propositions}\label{section:notation}
Our analysis takes place on a real pre-Hilbert space $\HH$,  with an
inner product denoted by $\langle\cdot,\cdot\rangle$.  We denote by
$\DIm(\HH)$ the dimension of $\HH$, possibly  infinite, and we often
denote subspaces of $\HH$ by $\VV$ or $\WW$. The orthogonal
complement of $\VV$ is denoted by $\VV^{\perp}$. If $\VV$ is a
complete subspace of $\HH$ (in particular finite dimensional), then
we denote the orthogonal projection of $\HH$ onto $\VV$ by
$P_{\VV}$.
%
We denote the norm induced by the inner product on $\HH$ by
$\|\cdot\|$, and the distance between $x,y \in \HH$ by $\dist(x,y)$
or equivalently $\|x-y\|$. Similarly,
$\dist(x,\VV)=\|P_{\VV}(x)-x\|$ is the induced distance between
$x\in\HH$ and a complete subspace $\VV\subseteq \HH$.

We have chosen to work in the general setting of a pre-Hilbert space
in order to emphasize the independence of our current and subsequent
results~\cite{LW-part1,LW-part2} from the dimension of the ambient
space.

By $d$ we denote an intrinsic dimension of interest to us, where
$d\in\mathbb{N}$ and $d+1\leq \DIm(\HH)$. We also use the integer
$k\geq1$ according to our purposes. Whenever we use $d$ or $k$ and
do not specify their range, one can always assume that they are
positive integers and $k, d+1\leq \DIm(\HH)$.

If $f$ is defined on $\HH^k$, then we denote the evaluation of $f$
on the ordered set of vectors $v_1,\ldots,v_{k+1}\in\HH$ with $v_j$
removed by $f(v_1,\ldots,v_{j-1},v_{j+1},\ldots,v_{k+1})$. We remark
that we maintain this notation for all $1\leq j\leq k$, in
particular, $j=1$ and $j=k+1$.  Similarly, for $1\leq j\leq k$, then
$f(v_1,\ldots,v_{j-1},u,v_{j+1},\ldots,v_k)$ is $f$ evaluated on the
ordered set of $k$ vectors
$v_1,\ldots,v_{j-1},u,v_{j+1},\ldots,v_k\in\HH$, where $v_j$ is
replaced by $u$.  We may remove two vectors, $v_i$ and $v_j$, from
the ordered set $\{v_1,\ldots,v_{k+2}\}$ and denote the function $f$
evaluated on the resulting set by
$f(v_1,\ldots,v_{i-1},v_{i+1},\ldots,v_{j-1},v_{j+1},\ldots,v_{k+2})$,
regardless of the order of $i$ and $j$ and whether or not either is
$1$ or $k+2$.  In this case the convention is always that $i\not=j$.

For an arbitrary subset $K$ in $\HH$, we denote its diameter by
$\diam(K)$. If $\mu$ is a measure on $\HH$, we denote the support
of $\mu$ by $\Supp$.

We follow with specific definitions and corresponding propositions
according to topics.
\subsection{Special Subsets of $\mathbf{\HH}$}
For an affine subspace $\LL\subseteq \HH$, a point $x\in \LL$, and
an  angle $\theta$ such that $0\leq\theta\leq\pi/2$, we define the
{\em cone}, $\Cne(\theta,\LL,x)$, centered at $x$ on $\LL$ in the
following way
$$\Cne(\theta,\LL,x):=\{u\in \HH:
\dist(u,\LL)\leq \|u-x\| \cdot \sin(\theta)\}.$$ For an affine
subspace, $\LL\subseteq \HH$ and $h > 0$, we define the {\em tube}
of height $h$ on $\LL$, $\Tbe(\LL,h)$, as follows.
$$\Tbe(\LL,h):=\{u\in \HH : \dist(u,\LL)\leq
h\}.$$For $r>0$ and $x\in \HH$, we define the {\em ball} of radius
$r$ on $x$ to be
$$\BB(x,r):=\{u\in \HH : \|u-x\|\leq r\}.$$

\subsection{Sets Generated by Vectors}
If $v_1\ldots,v_k \in \HH$, then the {\em parallelotope} spanned by
these vectors is the set $$\Pllp
(v_1,\ldots,v_k):=\left\{\sum_{i=1}^k t_i v_i \,: 0\leq t_i \leq 1,
\ i=1,\ldots,k\right\}.$$  Similarly, the {\em polyhedral cone}
spanned by $v_1,\ldots,v_k$ has the form
$$\Cpoly
(v_1,\ldots,v_k):=\left\{\sum_{i=1}^k t_i v_i \,: t_i \geq 0, \
i=1,\ldots,k\right \}.$$
The {\em affine plane} through the vectors $v_1,\ldots,v_k$  is
defined by
$$\Affn (v_1,\ldots,v_k):=\left\{\sum_{i=1}^kt_i v_i \, : \,
\sum_{i=1}^k t_i=1,\  \, t_i\in \RR,\,i=1,\ldots,k \right\}.$$ The
{\em convex hull} of $v_1,\ldots,v_k$ is the set $$\Chull
(v_1,\ldots, v_k):=\, \Affn(v_1\ldots,v_k)\cap \Cpoly
(v_1,\ldots,v_k).$$  If $\SS$ is a finite subset of $\HH$, we denote
the span of $\SS$ by $\LS{\SS}$, and sometimes also by $\Sp(\SS)$.

\subsection{Determinants and Contents} \label{subsec:content}
If $\HH$ is finite-dimensional, $\DIm(\HH)=k$, and
$\Phi=\{\phi_1,\ldots,\phi_k\}$ is an arbitrary orthonormal basis
for $\HH$, then we denote by $\det\nolimits_{\Phi}$ the determinant
function with respect to $\Phi$, that is, the unique alternating
multilinear function such that
$\det\nolimits_{\Phi}(\phi_1,\ldots,\phi_k)=1$. The following
elementary property of the determinant will be fundamental in part
of our analysis and hence we distinguish it.
\begin{proposition} \label{prop:det}If $\DIm(\HH)=k$, $v_1,\ldots,v_k
\in\HH$ and $u\in \Affn(v_1,\ldots,v_k)$, then for any orthonormal
basis~$\Phi$
$$\det\nolimits_{\Phi}(v_1,\ldots,v_k)=\sum_{i=1}^k
\det\nolimits_{\Phi}(v_1,\ldots,v_{i-1},u,v_{i+1},\ldots,v_k).$$
\end{proposition}The arbitrary  choice of $\Phi$ will not
matter to us and thus will not be specified.  Indeed, our major
statements will involve only $|\det\nolimits_{\Phi}|$, or will be
related to Proposition~\ref{prop:det}, both of which are invariant
under any choice of orthonormal basis $\Phi$. For this reason we
will usually refer to ``the determinant'' and dispense with the
subscript $\Phi$, i.e., $\det \equiv \det\nolimits_{\Phi}$.

If $v_1,\ldots,v_k\in\HH$, we define the $k$-content of the
parallelotope $\Pllp(v_1,\ldots,v_k)$, denoted by
$M_k(v_1,\ldots,v_k)$, as follows:
\begin{equation}\label{equation:definition-content}M_k(v_1,\ldots,v_k)
:=\begin{cases}\displaystyle\det\nolimits_{\Phi}(v_1,\ldots,v_k),&\textup{if
}k=\DIm(\HH)\textup{ for fixed
}\Phi,\\\displaystyle\left[\det\left(\left\{\langle
v_i,v_j\rangle\right\}_{i,j=1}^k\right)\right]^\frac{1}{2},&\textup{if
}k<\DIm(\HH).\end{cases}\end{equation}
We note that if $k=\DIm(\HH)$, then the $k$-content may obtain
negative values, and that the absolute value of the $k$-content
can be expressed by the same formula for all $k \leq \DIm(\HH)$,
i.e.,
$$|M_k(v_1,\ldots,v_k)|=\left[\det\left(\left\{\langle
v_i,v_j\rangle\right\}_{i,j=1}^k\right)\right]^\frac{1}{2}.$$

\subsection{High-Dimensional Sine
Functions}\label{subsection:sine-definition} 
\begin{figure}[htbp]
\begin{center}
\includegraphics[width = 8cm,height=5cm]{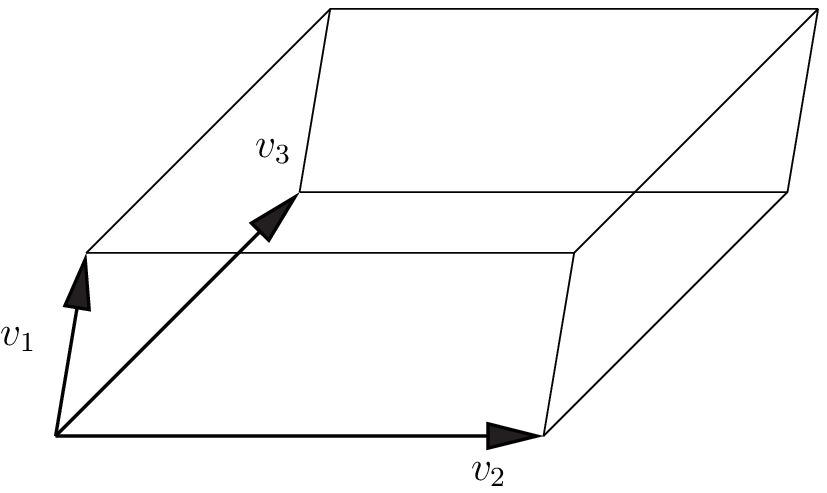}
\caption{\label{figure:def_sines}{\small Exemplifying the
computation of $\pds_0(v_1,v_2,v_3)$ and $\gds_0(v_1,v_2,v_3)$, when
$d=2$ and $\HH = \reals^3$: The figure shows the parallelepiped
spanned by $v_1$, $v_2$ and $v_3$. In this case,\newline
$\displaystyle \pds_0(v_1,v_2,v_3) = \frac{v_1 \bullet (v_2 \times
v_3)}{\|v_1\| \cdot \|v_2\| \cdot \|v_3\|}\ $ and
$\ \displaystyle \gds_0(v_1,v_2,v_3) = \frac{v_1 \bullet (v_2 \times
v_3)}{\sqrt{\|v_1 \times v_2 \| \cdot \|v_2 \times v_3 \| \cdot
\|v_1 \times v_3\|} }.$}}
\end{center}
\end{figure}

Using the definition of
$M_k$ in equation~(\ref{equation:definition-content}) and the
Euclidean norm on $\HH$, we define the functions
$\gds_0(v_1,\ldots,v_{d+1})$ and $\pds_0(v_1,\ldots,v_{d+1})$
respectively as

%
\begin{equation*}
\gds_0(v_1,\ldots, v_{d+1}) :=
\frac{M_{d+1}(v_1,\ldots v_{d+1})}
{\left(\prod^{d+1}_{j=1}M_{d}(v_1,\ldots v_{j-1},v_{j+1}\ldots
v_{d+1})\right)^{1/d}}
\end{equation*}
and
%
%
\begin{equation*}
\pds_0(v_1,\ldots,v_{d+1}) :=
\frac{M_{d+1}(v_1,\ldots,v_{d+1})}{\prod_{j=1}^{d+1}\|v_j\|}\,,%
%
\end{equation*}
where if either of the denominators above is zero (and thus the
numerator as well), then the corresponding function also obtains
the value zero.
We note that in the case of $d=1$, both functions are essentially
the ordinary sine functions. We exemplify this definition in
Figure~\ref{figure:def_sines}.

For these functions and their vector arguments
$v_1,\ldots,v_{d+1}$, we treat the point $0$ as a distinguished
vertex of the ($d+1$)-simplex through the vertices
$\{0,v_1,\ldots,v_{d+1}\}$. More generally, we may add a vertex
$w\in\HH$ other than $0$, and we define the functions
$\gds_w(v_1,\ldots,v_{d+1})$ and $\pds_w(v_1,\ldots,v_{d+1})$ for
vectors $v_1,\ldots,v_{d+1}$, $w\in\HH$ as follows:
\begin{equation} \label{eq:gsin_invariant}
\gds_w(v_1,\ldots,v_{d+1})=\gds_0(v_1-w,\ldots,v_{d+1}-w),
\end{equation}
and
\begin{equation} \label{eq:psin_invariant}
\pds_w(v_1,\ldots,v_{d+1})=\pds_0(v_1-w,\ldots,v_{d+1}-w).
\end{equation}
%

Whenever possible we refer to the functions $|\pds_w|$ and
$|\gds_w|$ so that we do not need to distinguish between the cases
$\DIm(\HH)=d+1$ and $\DIm(\HH) > d+1$. We mainly use the notation
$\pds_w$ or $\gds_w$ when $\DIm(\HH)=d+1$. In particular, we may use
the absolute values even if it is clear that $\DIm(\HH)>d+1$ and
thus the two sine functions are nonnegative.

We frequently use the following elementary property of $\pds_w$
and $\gds_w$, whose proof is included  in
Appendix~\ref{app:spherical}.
\begin{proposition} \label{prop:spherical} The functions
$|\pds_0|$ and $|\gds_0|$ defined on $\HH^{d+1}$ are invariant under
orthogonal transformations of $\HH$ and non-zero dilations of their
arguments. Moreover, if $\DIm(\HH)=d+1$, then $\pds_0$ and $\gds_0$
are invariant under dilations by positive coefficients.
%
\end{proposition}
%

%
%

Finally, we describe a generalized law of sines for $\gds$ following
Eriksson~\cite{eriksson} (see also~Barto\v{s}~\cite{bartos}):
\begin{proposition}
\label{prop:law} If $\{0,v_1,\ldots,v_{d+1}\}\subseteq \HH$ are
vertices of a non-degenerate $(d+1)$-simplex, then for all $1\leq
i\neq j\leq d+1$:
\begin{multline}\nonumber\frac{|\gds_0(v_1,\ldots,v_{d+1})|^d}{M_d(v_1-v_{d+1},\ldots,v_d-v_{d+1})}=\\
\frac{|\gds_{v_i}(v_1,\ldots,v_{i-1},0,v_{i+1},\ldots,v_{d+1})|^d}
{M_d(v_1-v_j,\ldots,v_{j-1}-v_j,v_{j+1}-v_j,\ldots,v_{i-1}-v_j,-v_j,v_{i+1}-v_j,\ldots,v_{d+1}-v_j)}\,.
\end{multline}
\end{proposition}
The proof follows from the definition of $|\gds_0|$. A reformulation
of this law is the invariance of the function
$|\gds_u(v_1,\ldots,v_{d+1})|/{M_d(v_1-u,\ldots,v_{d+1}-u)^{1/d}}$
with respect to permuting its arguments, $u$ included.

\subsection{ Elevation, Maximal Elevation, and Dihedral
Angles}\label{subsection:angle-definitions}

 For  a complete and non-trivial subspace $\WW\subseteq\HH$
and $u\in\HH\setminus\{0\}$, we define the {\em elevation angle} of
$u$ with respect to $\WW$ to be the smallest angle that $u$ makes
with any element $w\in\WW\setminus\{0\}$, and we denote this angle
by $\theta(u,\WW)$. More formally, in this case
$$\theta(u,\WW)=\min_{w\in\WW\setminus\{0\}}\left\{
\arccos\left(\left\langle \frac{u}{\|u\|}\,,
\frac{w}{\|w\|}\right\rangle\right)\right\}.$$If $u=0$, then we take
$\theta(0,W)=0$. We call the sines of these angles {\em elevation
sines} and note the following formula for computing them:
\begin{equation}\label{eq:elevation
sine}
\sin(\theta(u,\WW)) = \frac{\dist(u,\WW)}{\|u\|}\,.
\end{equation}

If $\VV$ is a complete subspace of $\HH$ and $v_1,v_2\in\HH$, we
define the {\em maximal elevation angle} of $v_1$ and $v_2$ with
respect to $\VV$, denoted by $\Theta(v_1,v_2,\VV)$, as
follows:\begin{equation}\Theta(v_1,v_2,\VV)=\max\{\theta(v_1,\VV),\theta(v_2,\VV)\}.\end{equation}

Given finite dimensional subspaces $\WW$ and $\VV$ of $\HH$ such
that $\DIm(\WW)=\DIm(\VV)$ and
$\DIm\left(\WW\cap\VV\right)=\DIm(\WW)-1$,  we define the {\em
dihedral angle} between $\WW$ and $\VV$ along $\WW\cap\VV$ to be the
acute angle between the normals of $\WW\cap\VV$ in $\WW$ and $\VV$.
We denote this angle by $\alpha(\WW,\VV)$. We call the sines of such
angles {\em dihedral sines} and note the following formula for
computing them:
\begin{equation}\label{eq:dihedral
sine}\sin(\alpha(\WW,\VV))=\frac{\dist(w,\VV)}{\dist(w,\WW\cap\VV)}
=\frac{\dist(v,\WW)}{\dist(v,\WW\cap\VV)},\ \textup{ for all }w\in
\WW\setminus\VV\textup{ and } v\in\VV\setminus\WW.\end{equation}

\subsection{Product Formulas for the High-Dimensional Sine Functions}

Two of the most useful properties of the high-dimensional sine
functions are their decompositions as products of lower-dimensional
sines. For $v_1,\ldots,v_{d+1}\in \HH$ and
$S=\{v_1,\ldots,v_{d+1}\}$, we formulate those decompositions as
follows.
\begin{proposition}\label{prop:prodg}$|\gds_0(v_1,\ldots,v_{d+1})|^d=
\left(\prod_{i=1}^d\sin\big(\AVU{\SS}{v_{d+1}}{v_i}\big)\right)\cdot
| \gdms_0(v_1,\ldots,v_d)|^{d-1}.$
\end{proposition}

\begin{proposition}\label{prop:prodp}$|\pds_0(v_1,\ldots,v_{d+1})|=\sin\big(\TV{\SS}{v_{d+1}}
{v_{d+1}}\big)\cdot |\pdms_0(v_1,\ldots,v_d)|.$\end{proposition}

Proposition~\ref{prop:prodg} was established
in~\cite[equation~7]{eriksson}, and Proposition~\ref{prop:prodp} can
be established given the fact that
\begin{eqnarray}\label{eq:contentformula}|M_{d+1}(v_1,\ldots,v_{d+1})| & = &
\dist(v_{d+1},\LV{\SS}{v_{d+1}})\cdot M_d(v_1,\ldots,v_d)\\
& = & \|v_{d+1}\|\cdot\sin\big(\TV{\SS}{v_{d+1}}{v_{d+1}}\big)\cdot
M_d(v_1,\ldots,v_d).\notag\end{eqnarray}


\section{Functional Identities for High-Dimensional Sine
Functions}\label{section:identities}

Throughout this section we assume that $\DIm(\HH)=d+1$ and
formulate identities for $\pds$ and $\gds$. We  denote the vectors
used for the arguments of the latter functions by $u$, $v_1,
\ldots , v_{d+1} \in \HH$, and assume the following:
$\{v_1,\ldots,v_{d+1}\}$ is a basis for $\HH$,
$u\in\Cpoly(v_1,\ldots,v_{d+1})$, and $u$ is not a scalar multiple
of any of the individual basis vectors $v_1,\ldots,v_{d+1}$, in
particular, $u \neq 0$.

\begin{figure}[htbp]
\begin{center}
\includegraphics[width = 9cm,height=9cm]{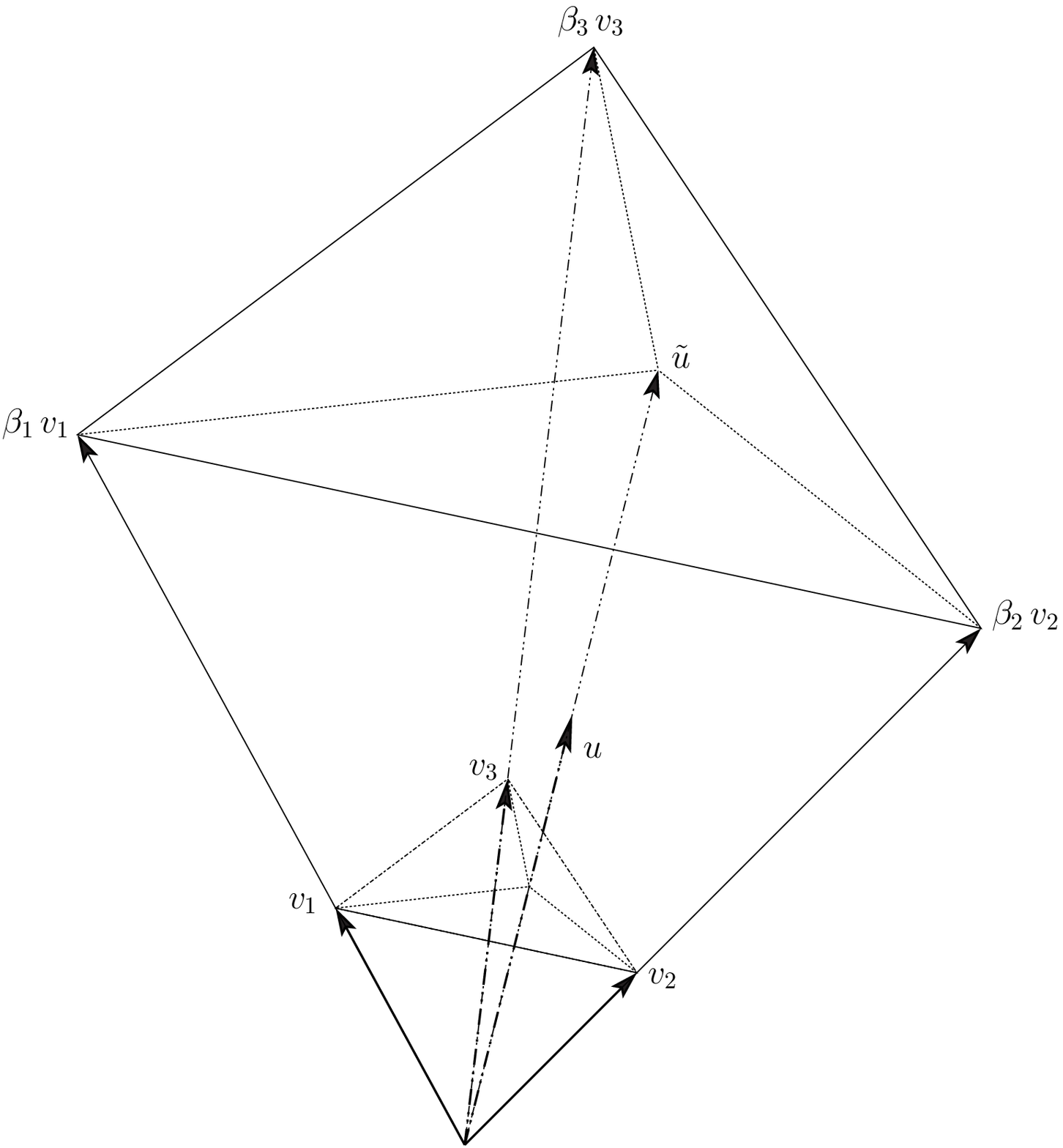}
\caption{\label{figure:identity_sines}{\small Exemplifying the
basic construction of this section, when $d=2$ and $\HH =
\reals^3$: We plot four particular vectors $v_1$, $v_2$, $v_3$ and
$u$ and note that in this special case $u$ is not contained in the
affine plane spanned by $v_1$, $v_2$ and $v_3$. We scale the
latter three vectors arbitrarily by the positive parameters
$\beta_1$, $\beta_2$ and $\beta_3$ and plot the resulting vectors.
We form $\tilde{u}$ by scaling $u$ so that it is in the affine
plane spanned by $\beta_1 v_1$, $\beta_2  v_2$ and $\beta_3
v_3$.}}
\end{center}
\end{figure}

%
The main elements of our identities are exemplified in
Figure~\ref{figure:identity_sines} and described as follows.
We introduce positive free parameters $\{\beta_i\}_{i=1}^{d+1}$, and
we note that $\Cpoly\left(\beta_1 v_1,\ldots,\beta_{d+1}
v_{d+1}\right)=\Cpoly(v_1,\ldots,v_{d+1})$. We express the vector
$u\in\Cpoly(v_1,\ldots,v_{d+1})$ as a linear combination of
$\{\beta_i v_i\}_{i=1}^{d+1}$ with coefficients
$\{\lambda_i\}_{i=1}^{d+1}$, that is,
\begin{equation}u=\sum_{i=1}^{d+1}\lambda_i\cdot
\beta_iv_i \,. \label{eq:lambda_i}
\end{equation}
We note that since $u\in\Cpoly(v_1,\ldots,v_{d+1})$ and $u \neq
0$, we have that $\sum_{i=1}^{d+1}\lambda_i
> 0$. We then define
%
%
\begin{equation}\label{eq:utilde}\tilde{u}:=
\left(\sum_{i=1}^{d+1}\lambda_i\right)^{-1} u,\end{equation}
and observe that
\begin{equation}\label{eq:u-in-affine}
\tilde{u}\in\Affn(\beta_1 v_1,\ldots,\beta_{d+1}
v_{d+1}).\end{equation} Finally, Proposition~\ref{prop:det} gives
the fundamental identity used to establish  all of the following
identities:
\begin{equation}\label{eq:detbeta}
\det(\beta_1v_1,\ldots,\beta_{d+1}v_{d+1})=\sum_{i=1}^{d+1}\det(\beta_1v_1,\ldots,\beta_{i-1}v_{i-1},
\tilde{u},\beta_{i+1}v_{i+1},\ldots,\beta_{d+1}v_{d+1}).
\end{equation}

In Subsection~\ref{subsection:identitypsine} we develop identities
for $\pds_0$ by direct application of the above equations.
Similarly, in Subsection~\ref{subsection:identitygsine} we develop
identities for $\gds_0$ following the same equations.  If $d=1$,
both identities for $\pds_0$ and $\gds_0$ reduce to a functional
equation satisfied by the sine function.  We characterize the
general Lebesgue measurable solutions of the corresponding equation
in Subsection~\ref{subsection:charmicheal}.

\subsection{Identities for $\mathbf{p_d
\text{sin}_0}$}\label{subsection:identitypsine}
Dividing both sides of equation~(\ref{eq:detbeta}) by
$\displaystyle\prod_{i=1}^{d+1}\|\beta_iv_i\|$, we obtain
$$\pds_0(\beta_1v_1,\ldots,\beta_{d+1}v_{d+1})=\sum_{i=1}^{d+1}P_i\cdot
\pds_0(\beta_1v_1,
\ldots,\beta_{i-1}v_{i-1},\tilde{u},\beta_{i+1}v_{i+1},\ldots,\beta_{d+1}v_{d+1}),
$$ where\begin{equation}\label{eq:polarcoefficients}
P_i \equiv
P_i\left(\{\beta_i\}_{i=1}^{d+1},\,\{v_i\}_{i=1}^{d+1},\,u\right)=\frac{\|\tilde{u}\|}{\|\beta_iv_i\|}.
\end{equation}
Applying either the law of sines or the formal definition of
$\pons$, we express the coefficients $P_i$ as
follows:\begin{equation}\label{eq:polarcoefficients2}
P_i
=\frac{\pons_0(-\beta_iv_i,
\tilde{u}-\beta_iv_i)}{\pons_0(\tilde{u},-\tilde{u}+\beta_iv_i)}\,.\end{equation}
%
 By the positive scale-invariance of
$\pds_0$ we obtain that
\begin{equation}\label{eq:polarequality}\pds_0(v_1,\ldots,v_{d+1})=\sum_{i=1}^{d+1}P_i\cdot
\pds_0(v_1,\ldots,v_{i-1},u,v_{i+1},\ldots,v_{d+1}).\end{equation}


By choosing different coefficients $\{\beta_i\}_{i=1}^{d+1}$ we can
obtain different identities for $\pds_0$. There are only $d$ degrees
of freedom in forming such identities due to the restriction of
equation~(\ref{eq:detbeta}).
In Subsection~\ref{subsec:simplexp} we will use the following choice
of $\{\beta_i\}_{i=1}^{d+1}$:
\begin{equation}\label{eq:def-beta-ratio-v_i}\beta_i=\frac{1}{\|v_i\|},\
i=1,\ldots,d+1.\end{equation}
The coefficients $\{P_i\}_{i=1}^{d+1}$, as described in
equation~(\ref{eq:polarcoefficients}), thus obtain the form,
\begin{equation}\label{eq:p_i-ratio-case}P_1=\ldots=P_{d+1}=\|\tilde{u}\|\,\end{equation}
and consequently equation~(\ref{eq:polarequality}) becomes
\begin{equation}\label{eq:psin-uniform}\pds_0(v_1,\ldots,v_{d+1})=
\|\tilde{u}\| \cdot
\sum_{i=1}^{d+1}\pds_0(v_1,\ldots,v_{i-1},u,v_{i+1},\ldots,v_{d+1}).\end{equation}

At last we exemplify the above identities when $d = 1$.
We denote the angle between $v_1$ and $u$ by $\alpha>0$, and the
angle between $u$ and $v_2$ by $\beta>0$, so that $\alpha+\beta$ is
the angle between $v_1$ and $v_2$. We note that by the two
assumptions of linear independence and $u \in \Cpoly (v_1,v_2)$ we
have that $\alpha+\beta<\pi$. We denote the angle between $-u$ and
$v_1-u$ by $\delta$, where $\beta<\delta<\pi-\alpha$. The parameter
$\delta$ represents the unique degree of freedom.

In this case, equations~(\ref{eq:polarcoefficients2})
and~(\ref{eq:polarequality}) reduce to the following trigonometric
identity:
%
%
%
\begin{equation}\label{eq:sin-identity}\sin(\alpha+\beta)=\frac{\sin(\alpha+\delta)}{\sin(\delta)}\cdot
\sin(\beta)\,+\,\frac{\sin(\delta-\beta)}{\sin(\delta)}\cdot
\sin(\alpha).\end{equation}
This identity generalizes to all $\alpha$, $\beta \in \reals$ and
$\delta \in \reals \setminus \pi \ints$. It was used in~\cite{MM05}
and is also very natural when establishing Ptolemy's theorem by
trigonometry.

Furthermore, equation~\eqref{eq:psin-uniform} reduces to the
trigonometric identity
$$\sin(\alpha+\beta)=\frac{\sin(\frac{\alpha+\beta}{2})}{\sin(\frac{\alpha-\beta}{2})}\cdot\big(\sin(\alpha)-\sin(\beta)\big),
$$
which can also be derived from equation~(\ref{eq:sin-identity}) by
setting $\delta=({\beta-\alpha})/{2}$.



\subsection{Identities for
$\mathbf{g_d\text{sin}_0}$}
\label{subsection:identitygsine} We now establish similar identities
for $\gds_0$.  Dividing both sides of equation~(\ref{eq:detbeta}) by
$\prod_{j=1}^{d+1}\left(M_d(\beta_1v_1,\ldots,\beta_{j-1}v_{j-1},\beta_{j+1}v_{j+1},
\ldots,\beta_{d+1}v_{d+1})\right)^{1/d}$ we obtain that
\begin{equation}\label{eq:1st-gsin-equality}\gds_0(v_1,\ldots,v_{d+1})=
\sum_{i=1}^{d+1}Q_i\cdot
\gds_0(\beta_1v_1,\ldots,\beta_{i-1}v_{i-1},\tilde{u},\beta_{i+1}v_{i+1},\ldots,\beta_{d+1}v_{d+1}),\end{equation}
where\begin{equation}\label{eq:orig-q_i} Q_i=
\left(\prod_{\substack{j=1\\j\not=i}}^{d+1}\frac{M_d(\beta_1v_1,\ldots,\beta_{j-1}v_{j-1},\beta_{j+1}v_{j+1},\ldots,\beta_{i-1}v_{i-1},
\tilde{u},\beta_{i+1}v_{i+1},\ldots,\beta_{d+1}v_{d+1})}{M_d(\beta_1v_1,\ldots,\beta_{j-1}v_{j-1},\beta_{j+1}v_{j+1},
\ldots,\beta_{d+1}v_{d+1})}\right)^{1/d}.\end{equation} By the
positive scale-invariance of $\gds_0$, we rewrite
equation~(\ref{eq:1st-gsin-equality}) as
\begin{equation}\label{eq:2nd-gsin-equality}\gds_0(v_1,\ldots,v_{d+1})=\sum_{i=1}^{d+1}Q_i\cdot
\gds_0(v_1,\ldots,v_{i-1},u,v_{i+1},\ldots,v_{d+1}).\end{equation}

We can express the coefficients $Q_i$ in  different ways.  First, we
note that
\begin{equation}\label{eq:1st-q_i}Q_i=\prod_{\substack{j=1\\j\not=i}}^{d+1}\frac{\gds_{\beta_iv_i}(\beta_1v_1,\ldots,\beta_{j-1}v_{j-1},0,\beta_{j+1}v_{j+1},\ldots,\beta_{i-1}v_{i-1},\tilde{u},
\beta_{i+1}v_{i+1},\ldots,\beta_{d+1}v_{d+1})}{\gds_{\tilde{u}}(\beta_1v_1,\ldots,\beta_{j-1}v_{j-1},0,\beta_{j+1}v_{j+1},\ldots,\beta_{d+1}v_{d+1})}.
\end{equation}  The fact that the absolute values of both
equations~(\ref{eq:orig-q_i}) and~(\ref{eq:1st-q_i}) are the same
follows from the generalized law of sines (see
Proposition~\ref{prop:law}). Moreover, the terms
$\{Q_i\}_{i=1}^{d+1}$ in equation~(\ref{eq:1st-q_i}) are positive
(see Appendix~\ref{app:sign_q_i}), as are the corresponding terms of
equation~(\ref{eq:orig-q_i}).

A different expression for $\{Q_i\}_{i=1}^{d+1}$ can be obtained
as follows. We set $\SS=\{v_1,\ldots,v_{d+1}\}$ and notice that
equation~(\ref{eq:contentformula}) implies that for all $1 \leq i
< j \leq d+1:$
$$\frac{M_d(\beta_1v_1,\ldots,\beta_{j-1}v_{j-1},\beta_{j+1}v_{j+1},\ldots,\beta_{i-1}v_{i-1},
\tilde{u},\beta_{i+1}v_{i+1},\ldots,\beta_{d+1}v_{d+1})}{M_d(\beta_1v_1,\ldots,\beta_{j-1}v_{j-1},\beta_{j+1}v_{j+1},
\ldots,\beta_{d+1}v_{d+1})}=\frac{\dist\left(\tilde{u},\LVU{\SS}{v_i}{v_j}\right)}{\dist\left(\beta_iv_i,\LVU{\SS}
{v_i}{v_j}\right)}\,.$$
Therefore, the coefficients $Q_i$,
$i=1,\ldots,d+1$, have the form
\begin{equation}\label{eq:4th-q_i}Q_i=\prod_{\substack{j=1\\j\not=i}}^{d+1}\left(\frac{\dist\left(\tilde{u},\LVU{\SS}{v_i}{v_j}\right)}{\dist\left(\beta_iv_i,\LVU{\SS}
{v_i}{v_j}\right)}\right)^{1/d}.\end{equation}
By further application of equation~(\ref{eq:elevation sine}), we
obtain that
\begin{equation}\label{eq:5ath-q_i}Q_i=\frac{\|\tilde{u}\|}{\|\beta_iv_i\|}\cdot
\prod_{\substack{j=1\\j\not=i}}^{d+1}\left(\frac{\sin\left(\TVU{\SS}{v_i}{v_j}{\tilde{u}}\right)}
{\sin\left(\TVU{\SS}{v_i}{v_j}{\beta_iv_i}\right)}\right)^{1/d}\,.
\end{equation}
It thus follows from equations~(\ref{eq:polarcoefficients})
and~(\ref{eq:5ath-q_i}) that
\begin{equation}\label{eq:5th-q_i}Q_i=P_i\cdot
\prod_{\substack{j=1\\j\not=i}}^{d+1}\left(\frac{\sin\left(\TVU{\SS}{v_i}{v_j}{\tilde{u}}\right)}
{\sin\left(\TVU{\SS}{v_i}{v_j}{\beta_iv_i}\right)}\right)^{1/d}\,.
\end{equation}

There are different possible choices for the parameters
$\{\beta_i\}_{i=1}^{d+1}$, and we present a specific choice and its
consequence in Subsection~\ref{subsubsection:simplexg}.
\subsection{Characterization of the Solutions of the One-Dimensional
Identity}\label{subsection:charmicheal} When $d=1$, the identities
of Subsections~\ref{subsection:identitypsine}
and~\ref{subsection:identitygsine} can be reduced to
equation~(\ref{eq:sin-identity}). That is, $f(x)=\sin(x)$ satisfies
the functional equation
\begin{equation}\label{eq:sin-k-identity}f(\alpha+\beta)=\frac{f(\alpha+\delta)}{f(\delta)}\cdot
f(\beta)\,+\,\frac{f(\delta-\beta)}{f(\delta)}\cdot f(\alpha) \ \
\text{ for all } \alpha,\beta \in \mathbb{R}, \ \delta \in \reals
\setminus f^{-1}(0)\,.
%
\end{equation}
%
%


We show here that the most general Lebesgue measurable solutions of
equation~(\ref{eq:sin-k-identity})
are multiples of the generalized sine functions on spaces of
constant curvature~\cite{polar_duality}, i.e., functions of the
form $c\cdot s_k(x)$, where
$$s_k(x)=\begin{cases} \frac{\sin(\sqrt{k}x)}{\sqrt{k}}, &
\textup{ if }k>0,\\ \quad x, & \textup{ if
} k=0,\\
\frac{\sinh(\sqrt{-k}x)}{\sqrt{-k}}, & \textup{ otherwise.
}\end{cases}$$

We remark that equation~(\ref{eq:sin-k-identity}) is almost
identical to an equation suggested by Mohlenkamp and
Monz{\'o}n~\cite[equation~(5)]{MM05}, but has a different set of
solutions. It is also closely related to Carmichael's
equation~\cite[Section~2.5.2, equation~(1)]{aczel} as the proof of
the following proposition shows.

We denote the set of multiples of generalized sine functions on
spaces of constant curvature by $\mathcal{S}$, that is,
\begin{equation}
\label{eq:sine_set} \nonumber {\mathcal{S}} = \big\{c\cdot s_k(x):\
c, k \in\RR\big\}.
\end{equation}
Using this notation, we formulate the main result of this
subsection:
\begin{theorem}\label{thm:charmicheal} The
set of all Lebesgue measurable functions satisfying
equation~(\ref{eq:sin-k-identity}) coincides with $\mathcal{S}$.
\end{theorem}

\begin{proof}
Clearly the elements of $\mathcal{S}$ satisfy
equation~(\ref{eq:sin-k-identity}). We thus assume that $f$ is a
Lebesgue measurable function satisfying
equation~(\ref{eq:sin-k-identity}) and show that $f \in
\mathcal{S}$. We denote the set of zeros of $f$ by $f^{-1}(0)$.
Since $f = 0$ is an element of $\mathcal{S}$ (obtained by setting
$c=0$ in equation~\eqref{eq:sine_set}), we also assume that $f \neq
0$ and in particular $\RR\setminus f^{-1}(0)$ is not empty.

We first observe that $f(0) = 0$. Indeed, by arbitrarily fixing
$\delta \in \RR\setminus f^{-1}(0)$ and setting $\alpha=-\beta$ in
equation~(\ref{eq:sin-k-identity}) we obtain that
\begin{equation}\label{eq:mol1}
f(0)=\frac{f(\delta-\beta)}{f(\delta)}\left(f(\beta)+f(-\beta)\right).\end{equation}
Setting also $\beta=0$, we get that $f(0)=0$.

Next, we show that the set $f^{-1}(0)$ has measure zero. We first
note that it is closed under addition. Indeed, if $\alpha$, $\beta
\in f^{-1}(0)$, then equation~(\ref{eq:sin-k-identity}) implies that
$f(\alpha + \beta) = 0$. Now, assuming that $f^{-1}(0)$ has positive
measure and applying a classical result of
Steinhaus~\cite[Theorem~6]{aczel-dhombres}, we obtain that
$f^{-1}(0)$, equivalently $f^{-1}(0)+f^{-1}(0)$, contains an open
interval. Then, if $0$ is an accumulation point of $f^{-1}(0)$, by
the additivity of zeros such an open interval extends to $\reals$,
that is, $f^{-1}(0) = \reals$, which is the case we excluded
($f=0$). Consequently, either $f^{-1}(0)$ has measure zero, or we
must have that $0$ is not an accumulation point of $f^{-1}(0)$ and
$f^{-1}(0)$ contains an open interval. However, this latter case
results in a contradiction as we show next.


%
Setting $\beta = \gamma \in\RR \setminus f^{-1}(0)$, $\delta \in
\RR \setminus f^{-1}(0)$, and $\alpha = \lambda \in f^{-1}(0)$ in
equation~(\ref{eq:sin-k-identity}), we get the formal relation:
\begin{equation}
\nonumber \frac{f(\gamma + \lambda)}{f(\gamma)} =  \frac{f(\delta +
\lambda)}{f(\delta)}\,.
\end{equation}
We thus conclude that for any $\lambda \in f^{-1}(0)$ there exists
a constant $C(\lambda) \in \reals\setminus\{0\}$ such that
\begin{equation}
\label{eq:gamma_delta} f(\gamma + \lambda) = C(\lambda) \cdot
f(\gamma) \ \text{ for all } \gamma \in \reals\,.
\end{equation}%
This equation implies that the case where $0$ is not an
accumulation point of $f^{-1}(0)$ and $f^{-1}(0)$ contains an open
interval cannot exist. We therefore conclude that $f^{-1}(0)$ has
measure zero.

Using the fact that $f^{-1}(0)$ is a null set and combining it
with equation~\eqref{eq:mol1}, we can show that $f$ is an odd
function. Indeed, if $f$ is not odd, then there exists a
$\beta\in\RR$ such that
$$f(\beta)+f(-\beta)\not=0\,,$$ and
by equation~(\ref{eq:mol1}) we have that $f(\delta-\beta)=0$ for
all $\delta\in\RR\setminus f^{-1}(0)$. Hence,
$$\RR\setminus f^{-1}(0)\,-\,\beta\subseteq f^{-1}(0)\,,$$ however this set inequality
contradicts the fact that $f^{-1}(0)$ is null. Therefore, $f$ is
odd.

%

We next observe that
\begin{equation}
\label{eq:state_f} \text{if } \lambda \in f^{-1}(0)\,, \text{ then }
|f(\gamma+\lambda)| = |f(\gamma)| \, \text{ for all }  \gamma \in
\reals \,.
\end{equation}
Indeed, fixing $\lambda \in f^{-1}(0)$ and replacing $\gamma$ with
$\gamma - \lambda$ in equation~\eqref{eq:gamma_delta}, we have
that
$$f(\gamma) = C(\lambda) \cdot f(\gamma-\lambda) \ \text{ for all } \gamma \in \reals \,.$$
Also, replacing $\gamma$ with $-\gamma$ in
equation~\eqref{eq:gamma_delta} and using the fact that $f$ is
odd, we obtain that
$$f(\gamma-\lambda) = C(\lambda) \cdot f(\gamma) \ \text{ for all } \gamma \in \reals \,.$$
The above two equations imply that $|C(\lambda)|=1$ for all
$\lambda\in f^{-1}(0)$ (the case $C(\lambda)=0$ is excluded by the
assumption that $f \neq 0$). Equation~\eqref{eq:state_f} then
follows from equation~\eqref{eq:gamma_delta}.

At last,
setting $\delta=\beta-\alpha$ in
equation~(\ref{eq:sin-k-identity}) and using the fact that $f$ is
odd, we get that
\begin{equation} \label{eq:pre_main} f(\alpha+\beta)=\frac{f(\beta)}{f(\beta-\alpha)}\cdot
f(\beta)-\frac{f(\alpha)}{f(\beta-\alpha)}\cdot f(\alpha) \ \text{
for all } \alpha, \, \beta \in \reals \text{ such that }
\beta-\alpha \in \reals \setminus f^{-1}(0)\,.
\end{equation}
Moreover, setting $\lambda=\alpha-\beta$ and $\gamma = \beta$ in
equation~\eqref{eq:state_f}, we obtain that
\begin{equation}
\label{eq:statement_f} |f(\alpha)| = |f(\beta)|
\ \text{ for all } \alpha, \, \beta \in \reals \text{ such that }
\alpha - \beta \in  f^{-1}(0)\,.
\end{equation}
Equations~\eqref{eq:pre_main} and~\eqref{eq:statement_f}
imply that $f$ satisfies
Carmichael's equation, i.e.,
$$f(\alpha+\beta)\cdot f(\beta-\alpha)=f(\beta)^2-f(\alpha)^2.$$
Since the Lebesgue measurable solutions of this equation are the
elements of $\mathcal{S}$
(see e.g., \cite[Corollary~15]{aczel-dhombres}), we conclude that
$f \in \mathcal{S}$.
\end{proof}
\section{Simplex Inequalities for High-Dimensional Sine
Functions}\label{section:inequalities}

In this section we prove Theorem~\ref{theorem:relaxedboth}, that is,
we show that the functions $|\gds_0|$ and $|\pds_0|$ are
$d$-semimetrics. We establish it separately for each of the
functions in the following theorems.

\begin{theorem}\label{theorem:simplexg} If
$v_1,\ldots,v_{d+1}\in \HH$ and $u\in\HH\setminus\{0\}$, then
\begin{equation*}|\gds_0(v_1,\ldots,v_{d+1})|\leq
\sum_{i=1}^{d+1}|\gds_0(v_1,\ldots,v_{i-1},u,v_{i+1},\ldots,v_{d+1})|.\end{equation*}
\end{theorem}
\begin{theorem}\label{theorem:simplexp}If
$v_1,\ldots,v_{d+1}\in \HH$ and $u\in\HH\setminus\{0\}$, then
\begin{equation*}|\pds_0(v_1,\ldots,v_{d+1})|\leq
\sum_{i=1}^{d+1}|\pds_0(v_1,\ldots,v_{i-1},u,v_{i+1},\ldots,v_{d+1})|.\end{equation*}
\end{theorem}

The proofs of both theorems are parallel.  We first prove them
when $\DIm(\HH)=d+1$ by applying the identities developed in
Section~\ref{section:identities}.  We then notice two phenomena of
dimensionality reduction. The first is that projection reduces the
values of $|\pds_0|$ and $|\gds_0|$. The second is that if $u \in
\left(\Sp(\{v_1,\ldots,v_{d+1}\})\right)^\perp$, then the
corresponding simplex inequality for $|\pds_0|$ and $|\gds_0|$
reduces to a relaxed simplex inequality of one term and constant
1.
We
remark that the second phenomenon of dimensionality reduction is
not fully necessary for concluding the theorems, i.e., using the
regular simplex inequality is fine, but we find it worth
mentioning.

We prove Theorem~\ref{theorem:simplexg} in
Subsection~\ref{subsec:simplexg} and Theorem~\ref{theorem:simplexp}
in Subsection~\ref{subsec:simplexp}.

\subsection{The Proof of Theorem~\ref{theorem:simplexg}}
\label{subsec:simplexg}
\subsubsection{\textbf{The Case of
$\text{\bf dim}{\mathbf{(\HH)=d+1}}$}}\label{subsubsection:simplexg}
We establish the following proposition.
\begin{lemma}\label{lemma:simplexg} If $\DIm(\HH)=d+1$,
$\{v_1,\ldots,v_{d+1}\}\subseteq \HH$ and $u\in\HH\setminus\{0\}$,
then
\begin{equation}\label{eq:gsinesimplex} |\gds_0(v_1,\ldots,v_{d+1})|\leq
\sum_{i=1}^{d+1}|\gds_0(v_1,\ldots,v_{i-1},u,v_{i+1},\ldots,v_{d+1})|.
\end{equation}
\end{lemma}

\begin{proof}[Proof of Lemma \ref{lemma:simplexg}]
Let $\SS=\{v_1,\ldots,v_{d+1}\}$. If $\SS$ is linearly dependent,
then $|\gds_0(v_1,\ldots,v_{d+1})|=0$, and the inequality holds.
Similarly, if $u$ is scalar multiple of any of the individual basis
vectors $v_1,\ldots,v_{d+1}$, then the inequality holds as an
equality. Thus, we may assume that $\Sp(\SS) =\HH$ and that $u$ is
not a scalar multiple of any of the individual basis vectors
$v_1,\ldots,v_{d+1}$.

Furthermore, we may assume  that $u\in \Cpoly(v_1,\ldots,v_{d+1})$.
Indeed, if this is not the case, then we may apply the following
procedure. We express $u$ as a linear combination of the vectors
$\{v_i\}_{i=1}^{d+1}$ using the coefficients
$\{\lambda_i\}_{i=1}^{d+1}$:

\begin{equation*}
u=\sum_{i=1}^{d+1}\lambda_i\,
v_i=\sum_{i=1}^{d+1}|\lambda_i|\,\sign(\lambda_i)\, v_i, \textup{
where }\ \sum_{i=1}^{d+1}|\lambda_i|\not=0.
\end{equation*}
For all $1\leq i \leq d+1$, we let
$$\hat
v_i=\begin{cases}\sign(\lambda_i)\, v_i, & \textup{ if }
\lambda_i\not=0,\\\qquad v_i, & \textup{ otherwise.}\end{cases}$$
We note that $u=\sum_{i=1}^{d+1}| \lambda_i|\cdot\hat v_i$, and
therefore $u\in \Cpoly(\hat v_1,\ldots,\hat v_{d+1})$.  Moreover, by
the scale-invariance of the function $|\gds_0|$ we obtain that the
required inequality (equation~(\ref{eq:gsinesimplex})) holds if and
only if
$$|\gds_0(\hat v_1,\ldots,\hat v_{d+1})|\leq
\sum_{i=1}^{d+1}|\gds_0(\hat v_1,\ldots,\hat v_{i-1},u,\hat
v_{i+1},\ldots,\hat v_{d+1})|.$$
Thus it is sufficient to consider the case where $u\in
\Cpoly(v_1,\ldots,v_{d+1})$. We observe that this assumption and
equation~(\ref{eq:u-in-affine}) imply that
\begin{equation}
\label{eq:u_til_in_cvx_hull} \tilde{u} \in
\Chull(v_1,\ldots,v_{d+1})\,.
\end{equation}

We next obtain the  desired inequality by using
equation~(\ref{eq:2nd-gsin-equality}) together with the form of
$\{Q_i\}_{i=1}^{d+1}$ set in equation~(\ref{eq:4th-q_i}).  The
question is how to choose the positive coefficients
$\{\beta_i\}_{i=1}^{d+1}$ such that $Q_i\leq 1$, $i=1,\ldots,d+1$.
Avoiding a messy optimization argument, we will show that there is a
natural geometric choice for the parameters
$\{\beta_i\}_{i=1}^{d+1}$. Indeed, letting
$$\beta_i=M_{d}(v_1,\ldots, v_{i-1},v_{i+1},\ldots v_{d+1}),\ i=1,\ldots,d+1,$$
we have that for all $1\leq i\leq d+1$
%
\begin{multline}\label{eq:faces}
M_d(\beta_1 v_1,\ldots\beta_{i-1} v_{i-1},\beta_{i+1}
v_{i+1},\ldots,\beta_{d+1}v_{d+1})  =\\
M_d( v_1,\ldots, v_{i-1}, v_{i+1},\ldots, v_{d+1}) \cdot
\prod^{d+1}_{\substack{j=1\\j\not= i}}\beta_j
 =  \prod^{d+1}_{j=1}\beta_j =
%
\prod^{d+1}_{j=1} M_d(v_1,\ldots,v_{j-1},v_{j+1},\ldots,v_{d+1}).
\end{multline}
In particular, for the simplex with vertices
$\{0,\beta_1v_1,\ldots,\beta_{d+1}v_{d+1}\}$, we obtain equal
contents for all $d$-faces containing the vertex $0$.

Another geometric property of the resulting simplex is that if
$1\leq k\not=i\leq d+1$, then $\beta_kv_k$ and $\beta_iv_i$ are of
equal distance from the $(d-1)$-plane $\LVU{\SS}{v_i}{v_k}$. That
is,
$$\dist(\beta_k v_k,\LVU{\SS}{v_i}{ v_k})
=\dist(\beta_ i v_i,\LVU{\SS}{v_i}{v_k}),\textup{ where } 1\leq
k\not=i\leq d+1.$$
This is a direct result of equation~(\ref{eq:contentformula}) and
the fact that the $d$-dimensional contents of the relevant faces are
equal (recall equation~(\ref{eq:faces})). Then, denoting the common
distance for both $\beta_k v_k$ and $\beta_iv_i$ from
$\LVU{\SS}{v_i}{ v_k}$ by $d_{ik}$, we note that
$$\{\beta_1 v_1,\ldots,\beta_{d+1} v_{d+1}\} \subseteq
\Tbe\left(\LVU{\SS}{v_i}{v_k},d_{ik}\right) \ \text{ for all } \
1\leq k\not=i\leq d+1\,.$$
Since $\Tbe\left(\LVU{\SS}{v_i}{v_k},d_{ik}\right)$ is convex,
$$\Chull(\beta_1 v_1,\ldots,\beta_{d+1} v_{d+1})\subseteq
\Tbe\left(\LVU{\SS}{v_i}{v_k},d_{ik}\right) \ \text{ for all } \
1\leq k\not=i\leq d+1\,.$$
This observation together with equation~(\ref{eq:u_til_in_cvx_hull})
imply that
\begin{equation*}
\tilde{u}\in\Tbe\left(\LVU{\SS}{v_i}{v_k},d_{ik}\right) \ \textup{
for all } \ 1\leq k\not=i\leq d+1,
\end{equation*}
that is,
\begin{equation}
\label{eq:frac_dist_gsin}
\frac{\dist\left(\tilde{u},\LVU{\SS}{v_i}{v_k}\right)}{\dist\left(\beta_iv_i,\LVU{\SS}{v_i}{v_k}\right)}\leq
1 \quad \textup{for all } 1\leq k\not=i\leq d+1.
\end{equation}
It follows from equations~(\ref{eq:4th-q_i})
and~(\ref{eq:frac_dist_gsin}) that $0\leq Q_i\leq 1$ for all $1\leq
i\leq d+1$ and the desired inequality is concluded.
\end{proof}

\subsubsection{\textbf{Dimensionality Reduction I}}
We show that projections reduce the value of $|\gds_0|$.
\begin{lemma}\label{lemma:projg} If $\VV$ is a $(d+1)$-dimensional subspace of
$\HH$, $\{v_1,\ldots,v_d\}\subseteq\VV$, $u\in\HH$, and
$P_{\VV}:\HH\rightarrow \VV$ is the orthogonal projection onto
$\VV$, then
\begin{equation}\label{eq:projg}|\gds_0(v_1,\ldots,v_d,P_{\VV}(u))|\leq
| \gds_0(v_1,\ldots,v_d,u)|.\end{equation}
\end{lemma}
\begin{proof}[Proof of Lemma \ref{lemma:projg}]
We form the sets $\BB=\{v_1,\ldots,v_d\}$,
$\SS=\{v_1,\ldots,v_d,u\}$ and
$\tilde{\SS}=\{v_1,\ldots,v_d,P_{\VV}(u)\}$.
In order to conclude the lemma it is sufficient to prove the
following inequality for dihedral angles:
\begin{equation}\label{eq:dihed-sine-cont}
\sin\left(\AVU{\tilde{\SS}}{P_{\VV}(u)}{v_i}\right) \leq
\sin\left(\AVU{\SS}{u}{v_i}\right),\quad\textup{for all }1\leq i\leq
d.\end{equation}
Indeed, equation~(\ref{eq:projg}) is a direct consequence of both
equation~(\ref{eq:dihed-sine-cont}) and the product formula for
$|\gds_0|$ of Proposition~\ref{prop:prodg}.

In order to prove the bound of equation~(\ref{eq:dihed-sine-cont})
it will be convenient to use the following orthogonal projections,
while recalling that $\BB=\{v_1,\ldots,v_d\}$:
\begin{alignat*}{2}
P_{\BB} & :&& \ \HH \rightarrow \LS{\BB},\\
N_{\BB} & :&& \ \HH \rightarrow \left(\LS{\BB}\right)^{\perp}\cap\VV,\\
P_i & :&& \ \HH\rightarrow \LV{\BB}{v_i},\ 1\leq i\leq d,\\
N_i & :&& \ \HH\rightarrow (\LV{\BB}{v_i})^{\perp} \cap \LS{\BB},\
1\leq i\leq d.
\end{alignat*}
We also define
$$N_{\VV} := I - P_{\VV}.$$
We note that
$u=P_{\VV}(u)+N_{\VV}(u)=P_i(u)+N_i(u)+N_{\BB}(u)+N_{\VV}(u)$, for
all $1\leq i \leq d$.

If  $N_{\BB}(u) = 0$, then $P_{\VV}(u) = P_{\BB}(u)$ and the set
$\{v_1,\ldots, v_d,P_{\VV}(u)\}$ is linearly dependent.  Hence,
$|\gds_0(v_1,\ldots,v_d,P_{\VV}(u))| = 0$ and the inequality holds
in this case.

If  $N_{\BB}(u)\not= 0$, we apply equation~(\ref{eq:dihedral sine})
and  obtain that
\begin{equation}\label{eq:normalstuff1}\sin\left(\AVU{\tilde{\SS}}{P_{\VV}(u)}{v_i}\right)
=\frac{\dist\left(P_{\VV}(u),\LV{\tilde{\SS}}{P_{\VV}(u)}\right)}{\dist\left(P_{\VV}(u),\LVU{\tilde{\SS}}{P_{\VV}(u)}{v_i}\right)}=
\frac{\|N_{\BB}(u)\|}{\|N_{\BB}(u)+N_i(u)\|},\quad 1\leq i \leq
d\,,\end{equation}
and
\begin{equation}\label{eq:normalstuff2}\sin
\left(\AVU{\SS}{u}{v_i}\right) =
\frac{\dist\left(u,\LV{\SS}{u}\right)}{\dist\left(u,\LVU{\SS}{u}{v_i}\right)}
=
\frac{||N_{\BB}(u)+N_{\VV}(u)||}{||N_i(u)+N_{\BB}(u)+N_{\VV}(u)||},\quad
1\leq i \leq d\,.\end{equation}
For any fixed $1\leq i \leq d$, the vectors  $N_{\BB}(u),\
N_{\VV}(u)$, and $N_i(u)$ are mutually orthogonal, and therefore
\begin{equation} \label{eq:normalstuff3}\sin \left(\AVU{\SS}{u}{v_i}\right)=
\frac{\|N_{\BB}(u)\|}{\|N_i(u)+N_{\BB}(u)\|} \,
\sqrt{\frac{\frac{||N_{\VV}(u)||^2}{||N_{\BB}(u)||^2} +
1}{\frac{||N_{\VV}(u)||^2}{||N_{\BB}(u)||^2+||N_i(u)||^2} + 1}} \geq
\frac{||N_{\BB}(u)||}{||N_i(u)+N_{\BB}(u)||}\,.\end{equation}
Equation~(\ref{eq:dihed-sine-cont}) follows from
equations~(\ref{eq:normalstuff1}) and~(\ref{eq:normalstuff3}), and
thus the lemma is concluded.
\end{proof}
\subsubsection{\textbf{Dimensionality Reduction II}}
We show how to relax the simplex inequality stated in
Theorem~\ref{theorem:simplexg} in the following special case.
\begin{lemma}\label{lemma:projorthg}
If $\VV$ is a $(d+1)$-dimensional
subspace of $\HH$, $\{v_1,\ldots,v_{d+1}\}\subseteq\VV$, and
$u\in\VV^{\perp}\setminus\{0\}$,  then
\begin{equation}\label{eq:orthoboundgsin}|\gds_0(v_1,\ldots
 v_{d+1})|\leq|\gds_0(v_1,\ldots,v_{i-1},u,v_{i+1},\ldots,
 v_{d+1})|,\quad \textup{ for all } 1\leq i\leq d+1.\end{equation}
\end{lemma}
\begin{proof}[Proof of Lemma~\ref{lemma:projorthg}]
We assume without loss of generality that
$\Sp(\{v_1,\ldots,v_{d+1}\})=\VV$ (otherwise
equation~\eqref{eq:orthoboundgsin} follows trivially).
We define the sets
$S_i=\{v_1,\ldots,v_{i-1},u,v_{i+1},\ldots,v_{d+1}\}$ for all
$1\leq i\leq d+1$. Since $u\in\VV^{\perp}\setminus\{0\}$ we obtain
from equation~(\ref{eq:dihedral sine}) that
\begin{equation}\label{eq:sine-u-ortho}\sin\left(\AVU{\SS_i}{u}{v_j}\right)=1,\quad
\textup{for all } 1\leq j\not=i\leq d+1.\end{equation}
Combining equation~(\ref{eq:sine-u-ortho}) with the product formula
for $|\gds_0|$ (Proposition~\ref{prop:prodg}) we get the following
equality for all $1\leq i\leq d+1$,
\begin{multline}\label{eq:reduction-1}
|\gds_0(v_1,\ldots,v_{i-1},u,v_{i+1},\ldots,v_{d+1})|^d\\= |
\gdms_0(v_1,\ldots,v_{i-1},v_{i+1},\ldots,
v_{d+1})|^{d-1}\prod^{d+1}_{\substack{j=1\\j\not=i}} \sin
\left(\AVU{\SS_i}{u}{v_j}\right)\\= |
\gdms_0(v_1,\ldots,v_{i-1},v_{i+1},\ldots,
v_{d+1})|^{d-1}.\end{multline} By further application of the product
formula for $|\gds_0|$ we obtain that for all $1\leq i\leq d+1$,
\begin{equation}\label{eq:reduction-2}|\gds_0(v_1,\ldots v_{d+1})|^d  \leq
| \gdms_0(v_1,\ldots,v_{i-1},v_{i+1}\ldots,v_{d+1})|^{d-1}.
\end{equation}
Equation~\eqref{eq:orthoboundgsin} thus follows from
equations~(\ref{eq:reduction-1}) and~(\ref{eq:reduction-2}).
\end{proof}

\subsubsection{\textbf{Conclusion of
Theorem~\ref{theorem:simplexg}}}

Let $P$ denote the orthogonal projection from $\HH$ onto
$\Sp\{v_1,\ldots,v_{d+1}\}$. If $P(u) = 0$, then we conclude the
Theorem from Lemma~\ref{lemma:projorthg}.

If  $P(u)\ne 0$, then we conclude the theorem by applying
Lemmata~\ref{lemma:simplexg} and~\ref{lemma:projg} successively as
follows:
\begin{multline*} |\gds_0(v_1,\ldots,v_{d+1})| \leq
\sum^{d+1}_{i=1}|\gds_0(v_1,\ldots,
v_{i-1},P(u),v_{i+1},\ldots,v_{d+1} )|\\ \leq
\sum^{d+1}_{i=1}|\gds_0(v_1,\ldots
v_{i-1},u,v_{i+1}\ldots,v_{d+1})|. \qed\end{multline*}

\subsection{The Proof of Theorem~\ref{theorem:simplexp}}
\label{subsec:simplexp} Here we prove essentially the same three
lemmata of Subsection~\ref{subsec:simplexg} for the function
$|\pds_0|$.

\subsubsection{\textbf{The Case of $\text{\bf dim}\mathbf{(\HH)=d+1}$}}
We establish here the following proposition.
\begin{lemma}\label{lemma:simplexp}If $\DIm(\HH)=d+1$,
$v_1,\ldots,v_{d+1}\in\HH$ and $u\in\HH\setminus\{0\}$, then
\begin{equation} |\pds_0(v_1,\ldots,v_{d+1})|\leq
\sum_{i=1}^{d+1}|\pds_0(v_1,\ldots,v_{i-1},u,v_{i+1},\ldots,v_{d+1})|.
\end{equation}
\end{lemma}

\begin{proof}[Proof of Lemma~\ref{lemma:simplexp}]
Similarly as in the proof of Lemma~\ref{lemma:simplexg}, we can
assume that $v_1,\ldots,v_{d+1}$ are linearly independent, $u$ is
not a scalar multiple of any of the individual basis vectors
$v_1,\ldots,v_{d+1}$ and
%
$u\in\Cpoly(v_1,\ldots,v_{d+1})$. Using the choice of
$\{\beta_i\}_{i=1}^{d+1}$ specified in
equation~\eqref{eq:def-beta-ratio-v_i}, we have that $\| \beta_i v_i
\| \leq 1$ for all $1 \leq i \leq d+1$. In view of
equation~\eqref{eq:u-in-affine} we can extend this bound to
$\tilde{u}$, i.e., we have that $\| \tilde{u} \| \leq 1$. The lemma
then follows by combining equation~\eqref{eq:psin-uniform} with the
latter bound.
\end{proof}
\subsubsection{\textbf{Dimensionality Reduction I} }
We show that projections reduce the value of $|\pds_0|$.
\begin{lemma}\label{lemma:projp}If $\VV$ is a $(d+1)$-dimensional subspace of
$\HH$,
 $v_1,\ldots,v_d \in\VV$, $u\in\HH$, and
$P_{\VV}:\HH\rightarrow \VV$ is the orthogonal projection onto
$\VV$, then
\begin{equation}\label{eq:projp}|\pds_0(v_1,\ldots,v_d,P_{\VV}(u))|\leq
| \pds_0(v_1,\ldots,v_d,u)|.\end{equation}
\end{lemma}
\begin{proof}[Proof of Lemma \ref{lemma:projp}]We form the sets
$\BB=\{v_1,\ldots v_d\}\subseteq \VV$, $\SS=\{v_1,\ldots,v_d,u\}$
and $\tilde{\SS}=\{v_1,\ldots,v_d,P_{\VV}(u)\}$.
In order to conclude the lemma, it is sufficient to prove that
\begin{equation}\label{eq:control-proj-sine}
\sin\left(\TV{\tilde{\SS}}{P_{\VV}(u)}{P_{\VV}(u)}\right) \leq
\sin\left(\TV{\SS}{u}{u}\right).\end{equation} Indeed,
equation~(\ref{eq:projp}) is a direct consequence of
equation~(\ref{eq:control-proj-sine}) and the product formula for
$|\pds_0|$ (Proposition~\ref{prop:prodp}).

In order to prove equation~(\ref{eq:control-proj-sine}), it will be
convenient to use the following orthogonal projections:
\begin{alignat*}{2}
P_{\BB} & : && \ \HH\rightarrow \LS{\BB},\\
N_{\BB} & : && \ \HH\rightarrow (\LS{\BB})^{\perp} \cap \VV.
\end{alignat*}
We also define
$$
N_{\VV}  := I - P_{\VV}\,.\\
$$
We note that
$u=P_{\VV}(u)+N_{\VV}(u)=P_{\BB}(u)+N_{\BB}(u)+N_{\VV}(u)$.

If $N_{\BB}(u)=0$, then $P_{\VV}(u)=P_{\BB}(u)\in \LS{\BB}$, and
the inequality (equation~(\ref{eq:projp}))  holds trivially since
the set $\tilde{\SS}=\{v_1,\ldots,v_d,P_{\VV}(u)\}$ is linearly
dependent.

If $N_{\BB}(u)\not=0$, we apply equation~(\ref{eq:elevation sine})
to obtain that
\begin{equation*}\sin\left(\TV{\SS}{u}{u}\right)=\frac
{\dist\left(u,\LV{\SS}{u}\right)}{\|u\|}=\frac{\|N_{\BB}(u)+N_{\VV}(u)\|}
{\|P_{\BB}(u)+N_{\BB}(u)+N_{\VV}(u)\|}\,,\end{equation*}
and
\begin{equation*}\sin\left(\TV{\tilde{\SS}}{P_{\VV}(u)}{P_{\VV}(u)}\right)=\frac
{\dist\left(P_{\VV}(u),\LV{\tilde{\SS}}{P_{\VV}(u)}\right)}{\|P_{\VV}(u)\|}
=\frac{\|N_{\BB}(u)\|}{\|P_{\BB}(u)+N_{\BB}(u)\|}\,.\end{equation*}
Thus,
\begin{equation*}\sin\left(\TV{\SS}{u}{u}\right)=
\sin\left(\TV{\tilde{\SS}}{P_{\VV}(u)}{P_{\VV}(u)}\right)
\sqrt{\frac{1+\frac{\|N_{\VV}(u)\|^2}{\|N_{\BB}(u)\|^2}}
{1+\frac{\|N_{\VV}(u)\|^2}{\|P_{\BB}(u)+N_{\BB}(u)\|^2}}}\geq
\sin\left(\TV{\tilde{\SS}}{P_{\VV}(u)}{P_{\VV}(u)}\right).\end{equation*}
That is, equation~(\ref{eq:control-proj-sine}) is verified and the
lemma is concluded. \end{proof}
\subsubsection{\textbf{Dimensionality Reduction II} }
We show how to relax the simplex inequality stated in
Theorem~\ref{theorem:simplexp} in the following special case.
\begin{lemma}\label{lemma:projorthp}If $\VV$ is a $(d+1)$-dimensional
subspace of $\HH$, $v_1,\ldots,v_{d+1}\in\VV$, and
$u\in\VV^{\perp}\setminus\{0\}$, then
\begin{equation}\label{eq:projorthp}
|\pds_0(v_1,\ldots,v_{d+1})|\leq|
\pds_0(v_1,\ldots,v_{i-1},u,v_{i+1},\ldots,v_{d+1})|,
\quad\textup{for
}i=1,\ldots,d+1\,.
\end{equation}

\end{lemma}
\begin{proof}[Proof of Lemma~\ref{lemma:projorthp}]
We assume without loss of generality that
$\Sp(\{v_1,\ldots,v_{d+1}\}) = \VV$. We define the sets
$S_i=\{v_1,\ldots,v_{i-1},u,v_{i+1},\ldots,v_{d+1}\}$ for all $1\leq
i\leq d+1$. Since $u\in\VV^{\perp}\setminus\{0\}$ we obtain from
equation~(\ref{eq:elevation sine}) that
\begin{equation}\label{eq:sine-equality-1}\sin \left(\TV{\SS_i}{u}{u}\right) =
1 ,\quad \textup{for all } 1\leq i \leq d+1\,. \end{equation}
Combining equation~(\ref{eq:sine-equality-1}) with the product
formula for $|\pds_0|$ (Proposition~\ref{prop:prodp}), we get the
following equality for all $i=1,\ldots,d+1$
\begin{equation}\label{eq:label1}
| \pds_0(v_1,\ldots,v_{i-1},u,v_{i+1},\ldots, v_{d+1})|
=|\pdms_0\left(v_1,\ldots,v_{i-1},v_{i+1},\ldots
,v_{d+1}\right)|\,.\end{equation}
By further application of the product formula for $|\pds_0|$, we
obtain that for all $i=1,\ldots,d+1$:
\begin{equation}\label{eq:whatever}
|\pds_0(v_1,\ldots,v_{d+1})|\leq|\pdms_0(v_1,\ldots,v_{i-1},v_{i+1},\ldots,v_{d+1})|.\end{equation}
Combining equations~(\ref{eq:label1}) and~(\ref{eq:whatever}) we
conclude equation~(\ref{eq:projorthp}). \qedhere
\end{proof}

\subsubsection{\textbf{Conclusion  of
Theorem~\ref{theorem:simplexp}}}


Let $P$ denote the orthogonal projection of $\HH$ onto
$\Sp(\{v_1,\ldots,v_{d+1}\})$. If $P(u) = 0$, then we conclude the
theorem from Lemma~\ref{lemma:projorthp}.

If $P(u)\ne 0$, then applying Lemmata~\ref{lemma:simplexp}
and~\ref{lemma:projp} successively we obtain that
\begin{multline*} |\pds_0(v_1,\ldots,v_{d+1})|  \leq
\sum^{d+1}_{i=1}|\pds_0(v_1,\ldots,
v_{i-1},P(u),v_{i+1},\ldots,v_{d+1} )| \\\leq
\sum^{d+1}_{i=1}|\pds_0(v_1,\ldots
v_{i-1},u,v_{i+1}\ldots,v_{d+1})|\,. \qed\end{multline*}
%

\section{ Ahlfors Regular Measures and Concentration Inequalities for the Polar Sine}
\label{section:relaxed} In this section we prove
Theorem~\ref{theorem:relaxedar}.  As explained in the introduction,
we interpret this theorem as indicating that the polar sine
$|\pds_w|$ satisfies a relaxed simplex inequality of two terms with
``high Ahlfors probability at all scales and locations''. Both
scales and locations are given by balls centered on $\Supp$, and
probabilities are given by scaling the $\gamma$-dimensional Ahlfors
regular measures of such balls, where $d-1<\gamma\leq d$.

\subsection{Notation, Definitions and Elementary Propositions}

For convenience of our notation, we assume that $w=0$ and
$0\in\Supp$, and thus establish most of the propositions for
$\pds_0$.  They can be generalized for $\pds_w$ via
equation~\eqref{eq:psin_invariant}.

Throughout this section we extensively use the definitions and
notation for elevation, maximal elevation, and dihedral angels
formulated in Subsection~\ref{subsection:angle-definitions}.  We
often fix $\SS=\{v_1,\ldots,v_{d+1}\}\subseteq\HH$.  If $0\leq
\epsilon \leq 1$ and $1\leq i\leq d+1$, then we denote by
$\Cne^i(\epsilon)$, the cone
$$\Cne^i(\epsilon)=\Cne\left(\epsilon\cdot\TV{\SS}{v_i}{v_i},\LV{\SS}{v_i},0\right).$$
If $0\leq \epsilon\leq 1$ and $1\leq i<j\leq d+1$, then we denote by
$\Cne^{i,j}(\epsilon)$ the set
\begin{equation}\label{def:Cneij}\Cne^{i,j}(\epsilon)=\Cne^i(\epsilon)\cap\Cne^j(\epsilon)\,.\end{equation}
If $1\leq i<j\leq d+1$, then we denote by $ \Theta_{i,j} $ the
following maximal elevation
angle\begin{equation}\label{eq:theta-i-j-def}\Theta_{i,j}=\Theta(v_i,v_j,\LVU{\SS}{v_i}{v_j}).\end{equation}

Throughout the rest of the paper we fix a real parameter $\gamma\in
\RR,\ d-1< \gamma \leq d$ (the most natural choice is $\gamma=d$)
and assume that $\HH$ is equipped with a $\gamma$-dimensional
Ahlfors regular measure, which we define as follows:
\begin{definition}\label{def:ARmeasure}  A locally finite Borel measure $\mu$
on $\HH$ is a  $\gamma$-dimensional Ahlfors regular measure if there
exists a constant $C$ such that for all $x\in \Supp$ and $0<r\leq
\diam(\Supp)$,
$$C^{-1}\cdot r^{\gamma}\leq \mu(\BB(x,r))\leq
C\cdot r^{\gamma}.$$
\end{definition}
We denote the smallest constant $C$ for which the inequality above
holds by $C_{\mu}$.  We refer to it as the {\em regularity
constant} of $\mu$. 

The following proposition and its immediate corollary, will be
useful for us. We prove them in Appendix~\ref{app:control_msr}.
\begin{proposition}\label{prop:artube}
If $\gamma> 1$, $m \in \mathbb{N}$ such that $1\leq m<\gamma$, $\mu$
a $\gamma$-dimensional Ahlfors regular measure on $\HH$ with
regularity constant $C_{\mu}$,
 $0 \leq \epsilon \leq 1$, and $\LL$ an $m$-dimensional affine subspace of
$\HH$, then for all $x\in\Supp\cap\LL$ and $0<r\leq \diam(\Supp)$
\begin{equation}
\label{eq:artube} \mu(\Tbe(\LL,\epsilon \cdot r )\cap \BB(x,r))\
\leq\ 2^{m+\frac{3\gamma}{2}}\cdot C_{\mu}\cdot \epsilon^{\gamma-m}
\cdot r^{\gamma}.\end{equation}
%
\end{proposition}

\begin{corollary} \label{cor:arcone}
If $\gamma > 1$, $m \in \mathbb{N}$ such that $1\leq m<\gamma$,
$\mu$ a $\gamma$-dimensional Ahlfors regular measure on $\HH$ with
regularity constant $C_{\mu}$, $0 \leq \theta\leq\pi/2$, and L an
$m$-dimensional affine subspace of $\HH$, then for all
$x\in\Supp\cap\LL$ and $0<r\leq \diam(\Supp)$
\begin{equation*}
\mu(\Cne(\theta,\LL,x)\cap \BB(x,r))\ \leq\
2^{m+\frac{3\gamma}{2}}\cdot C_{\mu}\cdot \sin(\theta)^{\gamma-m}
\cdot r^{\gamma}\,.
\end{equation*}
%
\end{corollary}

We
will frequently use the following elementary inequalities for the
one-dimensional sine:
\begin{lemma}\label{lemma:sine-ineq-1}If $0\leq\theta\leq \frac{\pi}{2}$ and
$0\leq c\leq 1$, then
\begin{equation} \label{eq:sine-ineq-1}c\cdot \sin(\theta)\leq \sin(c\cdot
\theta)\end{equation} and
\begin{equation}\label{eq:sin-ineq-3}\sin(c\cdot\theta)\leq\frac{\pi}{2}\cdot
c\cdot\sin(\theta).\end{equation}
\end{lemma}
Both inequalities can be derived by noting that the function
${\sin(c \, \theta)}/{(c\,\sin(\theta))}$ is increasing in
$\theta$ and thus obtains its lower bound, $1$, as $\theta$
approaches $0$ and its maximum value, bounded by $\pi/{2}$, at
$\theta=\frac{\pi}{2}$.

\subsection{Relationship between Conic Regions and Relaxed Inequalities for
the Polar Sine}
We establish here the following relation between the set
$U_C(\SS,0)$ defined in equation~(\ref{def:U-C})
and the intersection of various cones.

\begin{proposition}\label{prop:subsetu-C}If
$\SS=\{v_1,\ldots,v_{d+1}\}\subseteq\HH$, $C \geq 1$, $U_C(\SS,0)$
is the set defined in equation~(\ref{def:U-C}) with $w=0$, and
$\Cne^{i,j}(C^{-1})$ for $1\leq i<j\leq d+1$ are the intersections
of cones defined in equation~(\ref{def:Cneij}) with
$\epsilon=C^{-1}$, then
$$\HH\setminus\left(\bigcup_{1\leq i<j\leq
d+1}\Cne^{i,j}(C^{-1})\right)\subseteq
U_C(\SS,0).$$\end{proposition}

\begin{proof}  We note that $$\HH\setminus\bigcup_{1\leq i<j\leq
d+1}\Cne^{i,j}(C^{-1})=\bigcap_{1\leq i<j\leq
d+1}\left(\HH\setminus\Cne^{i,j}(C^{-1})\right),$$ and
$$U_C(\SS,0)=\bigcap_{1\leq i<j\leq d+1}U_C^{i,j}(\SS,0)\,,$$
where\begin{multline}\label{def:U-C-ij}U_C^{i,j}(\SS,0)=\Big\{u\in\HH:|\pds_0(v_1,\ldots,v_{d+1})|\leq
C\cdot\Big(|\pds_0(v_1,\ldots,v_{i-1},u,v_{i+1},\ldots,v_{d+1})|\
+\\|\pds_0(v_1,\ldots,v_{j-1},u,v_{j+1},\ldots,v_{d+1})|\Big)\Big\}.\end{multline}Therefore,
in order to conclude the proposition we only need to prove that
\begin{equation}\label{eq:contain-cones-U-C}\HH\setminus\Cne^{i,j}(C^{-1})\subseteq
U_C^{i,j}(\SS,0),\quad\textup{for all }1\leq i<j\leq
d+1.\end{equation}

If $u\in\HH\setminus\Cne^{i,j}(C^{-1})$ for some $i$ and $j$, where
$1\leq i<j \leq d+1,\ $ then either $u\in\HH\setminus\Cne^i(C^{-1})$
or $u\in\HH\setminus\Cne^j(C^{-1})$.  Assume without loss of
generality that $u\in\HH\setminus\Cne^i(C^{-1})$, then
\begin{equation}\label{eq:sine-bound-C-inv}\sin\left(\TV{\SS}{v_i}{u}\right)\geq\sin\left(C^{-1}\cdot\TV{\SS}{v_i}{v_i}\right).\end{equation}
Combining the product formula for $|\pds_0|$
(Proposition~\ref{prop:prodp}) with
equations~(\ref{eq:sine-bound-C-inv}) and~(\ref{eq:sine-ineq-1}), we
obtain that
$$C\cdot|\pds_0(v_1,\ldots,v_{i-1},u,v_{i+1},\ldots,v_{d+1})|\geq
| \pds_0(v_1,\ldots,v_{d+1})|.$$In particular, $u\in
U_C^{i,j}(\SS,0)$, and equation~(\ref{eq:contain-cones-U-C}), and
consequently the proposition, is concluded.
\end{proof}

\subsection{Controlling the Intersection of Two Cones}
The main part of the proof of Theorem~\ref{theorem:relaxedar} is to
show that for any $0\leq \epsilon\leq 1$ we can control the measure
of the sets $\Cne^{i,j}(\epsilon)$ which are defined in
equation~(\ref{def:Cneij}). We accomplish this by showing that such
sets are contained in specific cones on $(d-1)$-dimensional
subspaces of $\HH$, and then applying Corollary~\ref{cor:arcone} to
control the measure of the latter cones.  The crucial proposition is
the following.
\begin{proposition}\label{prop:v_1-v_2-spaces}If $0\leq s \leq 1,\ k\geq 3,\
\VV$ is a $k$-dimensional subspace of $\HH$,  $\VV_1$ and $\VV_2$
are two different $(k-1)$-dimensional subspaces of $\VV$,
$v_1\in\VV_1\setminus\VV_2$, $v_2\in\VV_2\setminus\VV_1$,  and
$\theta(v_1,\VV_2)$, $\theta(v_2,\VV_1)$, as well as
$\Theta(v_1,v_2,\VV_1\cap\VV_2)$ are the corresponding elevation and
maximal elevation angles, then
$$\Cne\left(\frac{2 \, s}{\sqrt{5} \, \pi}\cdot\theta\left(v_1,\VV_2\right),
\VV_2,0\right)\bigcap\Cne\left(\frac{2 \, s}{\sqrt{5} \,
\pi}\cdot\theta\left(v_2,\VV_1\right),\VV_1,0\right)
\subseteq\Cne\left(s\cdot
\Theta(v_1,v_2,\VV_1\cap\VV_2),\VV_1\cap\VV_2,0\right).$$
\end{proposition}

\begin{proof}  Note that
$\DIm\left(\VV_1\cap\VV_2\right)=k-2$, and that
$$\VV=\Sp\{v_1\}+\VV_2=\Sp\{v_2\}+\VV_1=\Sp\{v_1,v_2\}+\VV_1\cap\VV_2.$$
However, the above sum is not direct (the subspaces in the sum are
not mutually orthogonal).  We thus create a few orthogonal subspaces
which are expressed via the following orthogonal projections:
\begin{alignat*}{2}P_i &:&& \ \HH\rightarrow\VV_i,\  i=1,2,\\
N_i &:&& \ \HH\rightarrow\VV_i\cap (\VV_1\cap\VV_2)^{\perp},\ i=1,2,\\
P_{1,2} &:&& \ \HH \rightarrow\VV_1\cap\VV_2.\end{alignat*}
Note that\begin{equation}\label{eq:decompose-p}P_i=P_{1,2}+N_i,\quad
i=1,2,\end{equation} and consequently
\begin{equation}\label{eq:relation-pn}\left(I-P_2\right)\cdot
P_1=\left(I-P_2\right)\cdot N_1.\end{equation}

We denote
\begin{align*}
&\Cnetilde^1=\Cne\left(\frac{2 \, s}{\sqrt{5} \, \pi}\cdot\theta\left(v_2,\VV_1\right),\VV_1,0\right),\\
&\Cnetilde^2=\Cne\left(\frac{2 \, s}{\sqrt{5} \, \pi}\cdot\theta\left(v_1,\VV_2\right),\VV_2,0\right),\\
&\Cnetilde^{1,2}=\Cnetilde^1\cap\Cnetilde^2,\\
%
%
&\Ttilde_{1,2}=\Theta(v_1,v_2,\VV_1\cap\VV_2). \end{align*}
Following our notation and the definition of a cone, we need to
prove that
$$\|u-P_{1,2}(u)\|\leq
\sin\left(s\cdot\Ttilde_{1,2}\right)\cdot\|u\|, \ \text{ for all }
u\in\Cnetilde^{1,2},$$
or equivalently (via equation~(\ref{eq:decompose-p})),
%
\begin{equation}\label{eq:toconclude-np}\|N_1(u)+u-P_1(u)\|\leq
\sin\left(s\cdot\Ttilde_{1,2}\right) \cdot \|u\|, \ \text{ for all }
u\in\Cnetilde^{1,2}.\end{equation}
For $u\in\Cnetilde^{1,2}$, we will bound $\|N_1(u)\|$ and
$\|u-P_1(u)\|$ separately and then combine the two estimates to
conclude the above inequality and the current proposition.

Our bound for $\|u-P_1(u)\|$ is straightforward.  Indeed, if
$u\in\Cnetilde^{1,2}$ respectively, then $u\in\Cnetilde^1$, and by
the definition of $\Cnetilde^1$ as well as the application of
equation~(\ref{eq:sin-ineq-3}) we obtain that
\begin{equation}\label{eq:bound-u-minus-p-1}\|u-P_1(u)\|\leq\sin\left(\frac{2 \, s}{\sqrt{5} \, \pi}\cdot\theta(v_2,\VV_1)\right)
\cdot\|u\|\leq\frac{s}{\sqrt{5}}
\cdot\sin\left(\theta(v_2,\VV_1)\right)\cdot\|u\|.\end{equation}

Our bound for $\|N_1(u)\|$ has the following form:
\begin{equation}\label{eq:bound-n1-3}\|N_1(u)\|\leq\frac{s}{\sqrt{5}}\left[\sin(\theta(v_1,\VV_1\cap\VV_2))+
\sin(\theta(v_2,\VV_1\cap\VV_2))\right]\cdot\|u\|.\end{equation}
In order to verify it, we assume without loss of generality that
$N_1(u)\not=0$ and note that equation~(\ref{eq:dihedral sine})
implies the following relation
\begin{equation}\label{eq:dihedral-pn-1}\sin(\alpha(\VV_1,\VV_2))=\frac{\dist(N_1(u),\VV_2)}{\dist(N_1(u),\VV_1\cap
\VV_2)}=\frac{\|N_1(u)-P_2\cdot
N_1(u)\|}{\|N_1(u)\|}.\end{equation}Combining
equations~(\ref{eq:relation-pn}) and~(\ref{eq:dihedral-pn-1}) we
obtain that
\begin{equation}\label{eq:express-n-1}\|N_1(u)\|=\frac{\|P_1(u)-P_2\cdot
P_1(u)\|}{\sin\left(\alpha(\VV_1,\VV_2)\right)}.\end{equation}
We bound $\|P_1(u)-P_2\cdot P_1(u)\|$ as follows.
\begin{multline}\label{eq:ineq-p-n}
\|P_1(u)-P_2\cdot P_1(u)\|=\|(I-P_2)\cdot
P_1(u)\|\\
=\|(I-P_2)(u)-(I-P_2)\cdot (I-P_1)(u)\|
\leq\|u-P_2(u)\|+\|u-P_1(u)\|.
\end{multline}Equation~(\ref{eq:bound-u-minus-p-1}) gives a bound for
$\|u-P_1(u)\|$, and similarly we obtain that
\begin{equation}\label{eq:bound-u-minus-p-2}\|u-P_2(u)\|\leq
\frac{s}{\sqrt{5}}\cdot \sin(\theta(v_1,\VV_2)\cdot
\|u\|.\end{equation}Combining equations~(\ref{eq:bound-u-minus-p-1})
and~(\ref{eq:express-n-1})-(\ref{eq:bound-u-minus-p-2}) we get that
\begin{equation}\label{eq:bound-n1-2}\|N_1(u)\|\leq\frac{s}{\sqrt{5}}\left(\frac{\sin\left(\theta\left(v_1,\VV_2\right)\right)}{\sin(\alpha(\VV_1,
\VV_2))}+
\frac{\sin\left(\theta\left(v_2,\VV_1\right)\right)}{\sin(\alpha(\VV_1,\VV_2))}\right)\cdot\|u\|.\end{equation}
At last we note that equations~(\ref{eq:elevation sine})
and~(\ref{eq:dihedral sine}) imply
that$$\sin(\alpha(\VV_1,\VV_2))=\frac{\sin(\theta(v_1,\VV_2))}{\sin(\theta(v_1,\VV_1\cap\VV_2))}
=\frac{\sin(\theta(v_2,\VV_1))}{\sin(\theta(v_2,\VV_1\cap\VV_2))}.$$
Applying this identity in~(\ref{eq:bound-n1-2}), we achieve the
bound for $\|N_1(u)\|$ stated in equation~\eqref{eq:bound-n1-3}.

Finally, noting that $N_1(u)\perp (u-P_1(u))$ and applying the
bounds of equations~(\ref{eq:bound-u-minus-p-1})
and~(\ref{eq:bound-n1-3}) we obtain that
\begin{multline*}\|N_1(u)+u-P_1(u)\|^2=\|N_1(u)\|^2+\|u-P_1(u)\|^2\\\leq
\left(\frac{s}{\sqrt{5}}\right)^2\cdot 5\cdot
\sin^2\left(\Ttilde_{1,2}\right)\cdot\|u\|^2= s^2\cdot
\sin^2\left(\Ttilde_{1,2}\right)\cdot\|u\|^2
\leq\sin^2\left(s\cdot\Ttilde_{1,2}\right)\cdot\|u\|^2.\end{multline*}
Equation~(\ref{eq:toconclude-np}), and consequently the
proposition, is thus concluded.
\end{proof}

\begin{remark}
The proposition extends trivially to $k=2$, where the intersection
of two cones, centered around two vectors $w_1$ and $w_2$
respectively with opening angles less than half the angle between,
is the origin, which is a degenerate cone.
\end{remark}

Proposition~\ref{prop:v_1-v_2-spaces} implies the following
corollary:
\begin{corollary}\label{cor:cor-angles}
If $0\leq s\leq 1$, $2\leq k\leq \DIm(\HH)$, $1\leq i<j\leq d+1$,
$\SS=\{v_1,\ldots,v_k\}\subseteq\HH$ is a linearly independent set
and $\Theta_{i,j}$ as well as $\Cne^{i,j}\left(\frac{2\cdot
s}{\sqrt{5} \, \pi}\right)$ are defined by
equations~(\ref{eq:theta-i-j-def}) and~(\ref{def:Cneij})
respectively, then
\begin{equation}\Cne^{i,j}\left(\frac{2\cdot
s}{\sqrt{5} \,
\pi}\right)\subseteq\Cne\left(s\cdot\Theta_{i,j},\LVU{\SS}{v_i}{v_j},0\right).\end{equation}
\end{corollary}
Indeed, Corollary~\ref{cor:cor-angles} is obtained as a special case
of Proposition~\ref{prop:v_1-v_2-spaces} by setting $\VV=\LS{\SS}$,
$\VV_1=\LV{\SS}{v_i}$, and $\VV_2=\LV{\SS}{v_j}$, and noting that
$\VV_1\cap\VV_2=\LVU{\SS}{v_i}{v_j}$.

\subsection{Conclusion of Theorem~\ref{theorem:relaxedar}}

Theorem~\ref{theorem:relaxedar} follows directly from
Proposition~\ref{prop:subsetu-C} and Corollaries~\ref{cor:arcone}
and~\ref{cor:cor-angles}.

In view of equation~\eqref{eq:psin_invariant}, we note that it is
sufficient to prove the theorem when $w=0$ and $0 \in \Supp$. We
assume an arbitrary parameter $0<s\leq1$ and set
\begin{equation}\label{eq:set-C-zero}C=\frac{\sqrt{5} \, \pi}{2\cdot
s}\,.\end{equation}  At the end of the proof we further restrict
the values of $s$ from above and consequently restrict those of
$C$ from below.

Let $\SS=\{v_1,\ldots,v_{d+1}\}\subseteq\HH$, $\
0<r\leq\diam(\Supp)$, and $\Cne^{i,j}(C^{-1})$ be defined according
to equation~(\ref{def:Cneij}). We assume without loss of generality
that the set $\SS$ is linearly independent.
Proposition~\ref{prop:subsetu-C} implies that
$$\BB(0,r)\setminus\left(\bigcup_{1\leq i\not=j\leq
d+1}\Cne^{i,j}(C^{-1})\right)\subseteq \BB(0,r)\bigcap
U_{C}(\SS,0)\,.$$
Using the additivity and monotonicity of $\mu$, we
get
\begin{equation}\label{eq:u-C-measure}\mu\left(\BB(0,r)\cap
U_{C}(\SS,0)\right)\geq \mu\left(\BB(0,r)\right)\ -\ \sum_{1\leq
i<j\leq
d+1}\mu\left(\BB(0,r)\cap\Cne^{i,j}(C^{-1})\right).\end{equation}
Next, we combine Corollary~\ref{cor:cor-angles} together with
equation~(\ref{eq:set-C-zero}) to obtain that
\begin{equation}\label{eq:containingcone}\BB(0,r)\cap\Cne^{i,j}(C^{-1})
\subseteq
\BB(0,r)\cap\Cne\left(s\cdot\Theta_{i,j},\LVU{\SS}{v_i}{v_j},0\right)
\subseteq
\BB(0,r)\cap\Cne\left(s\cdot\frac{\pi}{2},\LVU{\SS}{v_i}{v_j},0\right).\end{equation}
Now, Corollary~\ref{cor:arcone}, Definition~\ref{def:ARmeasure}, and
equation~(\ref{eq:containingcone}) imply that for all $1\leq
i\not=j\leq d+1$,
\begin{multline}\label{eq:lastmeasure}\mu\left(\BB(0,r)\cap
\Cne^{i,j}(C^{-1})\right)\leq\mu\left(\BB(0,r)\cap\Cne\left(s\cdot\frac{\pi}{2},\LVU{\SS}{v_i}{v_j},0\right)\right)\\\leq
2^{\frac{3\,\gamma}{2}+d-1}\cdot
C_{\mu}^2\cdot\left(\sin\left(s\cdot\frac{\pi}{2}\right)\right)^{\gamma+1-d}
\cdot\mu(\BB(0,r))\,.
\end{multline}
Combining equations~(\ref{eq:u-C-measure})
and~(\ref{eq:lastmeasure}), we get that
\begin{equation}\frac{\mu\big(\BB(0,r)\cap
U_{C}(\SS,0)\big)}{\mu\left(\BB(0,r)\right)}\geq1\ -\
\binom{d+1}{2}\cdot 2^{\frac{3\,\gamma}{2}+d-1} \cdot C_{\mu}^2\cdot
\left(\sin\left(s\cdot\frac{\pi}{2}\right)\right)^{\gamma+1-d}.
\end{equation}

By setting the parameter $s$ so that
$$\binom{d+1}{2}\cdot 2^{\frac{3\,\gamma}{2}+d-1} \cdot C_{\mu}^2\cdot
\left(\sin\left(s\cdot\frac{\pi}{2}\right)\right)^{\gamma+1-d} \leq
\epsilon\,,$$
that is,
$$s \leq s_0^{\prime} = \frac{2}{\pi}\cdot\arcsin\left[\left(\frac{\epsilon}{2^{\frac{3\,\gamma}{2}+d-1}
\cdot
C_{\mu}^2\cdot\binom{d+1}{2}}\right)^{\frac{1}{\gamma+1-d}}\right]
,$$
we obtain that equation~\eqref{eq:theorem2} is satisfied for all $C
\geq C_0^\prime$, where
\begin{multline}
\label{eq:bound_c0} \nonumber C_0^\prime=\frac{\sqrt{5} \, \pi}{2\,
s_0^\prime}=\sqrt{5}\,\left(\frac{\pi}{2}\right)^2\,\left(\arcsin\left[
\left(\frac{\epsilon}{2^{\frac{3\,\gamma}{2}+d-1} \cdot
C_{\mu}^2\cdot\binom{d+1}{2}}\right)^{\frac{1}{\gamma+1-d}}\right]\right)^{-1}
\\
%
=
O\left(\left(2^{\frac{3\,\gamma}{2}+d} \cdot
C_{\mu}^2\cdot\binom{d+1}{2}\cdot
\epsilon^{-1}\right)^{\frac{1}{\gamma + 1 -d}}\right)
 \ \textup{ as } \epsilon
\rightarrow 0 \ \textup{ or } \ d \rightarrow \infty\,.
\end{multline}
The theorem is thus concluded, where $C_0^\prime$ provides an
upper bound for the best possible choice for the constant $C_0$.
%
\qed
\begin{remark}
\label{remark:gamma_large} Note that Theorem~\ref{theorem:relaxedar}
extends trivially to the case where $\gamma > d$. In fact, in this
case, it is possible to replace the set $U_C(\SS,w)$ by
\begin{multline*}U'_C(\SS,w)=\big\{u\in\HH:\
| \pds_w(v_1,\ldots,v_{d+1})|\leq \\
 C\cdot|\pds_w(v_1,\ldots,v_{i-1},u,v_{i+1},\ldots,v_{d+1})|,
\textup{ for all }1\leq i \leq d+1 \big\}.\end{multline*}
That is, if $\gamma > d$, then the polar sine satisfies a relaxed
simplex inequality of one term ``with high Ahlfors probability at
all scales and locations''. This fact
is a direct consequence of Corollary~\ref{cor:arcone} and analogues
of Proposition~\ref{prop:subsetu-C} and
equation~(\ref{eq:u-C-measure}) obtained by replacing $U_C(\SS,0)$
with $U'_C(\SS,0)$ and $\{\Cne^{i,j}(C^{-1})\}_{1\leq i < j \leq
d+1}$ with $\{\Cne^{i}(C^{-1})\}_{1\leq i < d+1}$.

Nevertheless, if $d-1<\gamma\leq d$, then one cannot replace the set
$U_C(\SS,w)$ in Theorem~\ref{theorem:relaxedar}  by
$U'_C(\SS,w)$.\end{remark}

\begin{remark} Let us slightly reformulate the above results so that they could be
directly applied in~\cite{LW-part1}. For
$\SS=\{v_1,\ldots,v_{d+1}\}$ as above, $C>0$, and an arbitrarily
fixed pair of indices $i$ and $j$, where $1\leq i<j\leq d+1$, we
form the set $U_{C}(\SS,i,j,0)$ as follows:
\begin{multline*}U_{C}(\SS,i,j,0)=\Big\{u\in\HH:\ |
\pds_0(v_1,\ldots,v_{d+1})|\leq
C\cdot\big(|\pds_0(v_1,\ldots,v_{i-1},u,v_{i+1},\ldots,v_{d+1})|\
+\\|\pds_0(v_1,\ldots,v_{j-1},u,v_{j+1},\ldots,v_{d+1})|\big)\Big\}.\end{multline*}
If $\gamma = d$ and $0< \epsilon < 1$, then for all $C \geq
C_0^{\prime \prime}$, where
$$
C_0^{\prime \prime} =
\sqrt{5}\cdot\left(\frac{\pi}{2}\right)^2\cdot\left(\arcsin
\left(\frac{\epsilon}{2^{\frac{5\,d}{2}-1} \cdot
C_{\mu}^2}\right)\right)^{-1}
,$$%
we have that
$$\frac{\mu\big(\BB(0,r)\cap
U_{C}(\SS,i,j,0)\big)}{\mu\left(\BB(0,r)\right)}\geq1 -
\epsilon.$$\end{remark}

\section{Conclusions and Further Directions}\label{section:final}
The work presented here
touches on both old and modern research. We would like to conclude
it by indicating various directions where one can extend it.

\subsection*{High-Dimensional Menger-type curvature.}
In~\cite{LW-part1,LW-part2} we build on the research presented here
to form $d$-dimensional Menger-type curvatures of any integer
dimension $d > 1$, and show how they characterize $d$-dimensional
uniform rectifiability of $d$-dimensional Ahlfors regular measures
on real separable Hilbert spaces.

\subsection*{The high-dimensional $\mathbf{A=B}$ Paradigm.}

Petkov{\v{s}}ek, Wilf, and Zeilberger~\cite{zeilberger} have
presented concrete strategies to prove various identities.
However, when considering the high-dimensional sine functions, it
is not clear whether a general mechanism exists. The product
formulas (Propositions~\ref{prop:prodg} and~\ref{prop:prodp})
simplify the representation of $\pds_0$ and $\gds_0$, but they do
not seem to provide sufficiently simple structure for
automatically proving general identities involving those
functions. We have demonstrated additional strategies for proving
identities of interest to us and inquire about other useful
identities and the strategies for proving them.

\subsection*{Solutions of high-dimensional
functional equations.} We have shown that
equation~(\ref{eq:sin-k-identity}) characterizes the generalized
sine function of spaces with constant one-dimensional curvature
among all Lebesgue measurable functions
(Theorem~\ref{thm:charmicheal}). It will be interesting to formulate
a theorem analogous to Theorem~\ref{thm:charmicheal} for the
high-dimensional functional equations described in
Section~\ref{section:identities}. In particular, we are interested
in the functional equation generalizing the combination of
equations~(\ref{eq:1st-gsin-equality}) and~(\ref{eq:1st-q_i}). We
could not identify any similar functional equation in the
substantial body of work on the subject (see e.g.,
\cite{aczel,aczel-dhombres} and references therein).
\subsection*{Other relaxed inequalities with high probability.}
%
We inquire about probabilistic
settings different than the one in here, where the polar sine satisfies a
relaxed simplex inequality of two terms, but not of one term, with high probability.
We also inquire about other probabilistic settings in which the
polar sine satisfies a relaxed simplex inequality of $p$ terms, $3
\leq p \leq d$, and not of $p-1$ terms, with high probability.
Moreover, we are curious about probabilistic settings where one
can obtain relaxed simplex inequalities for $|\gds|$ with high
probabilities.
\subsection*{Applications to data analysis.} Recently, researchers in
machine learning have been interested in multi-way clustering and $d$-way kernel
methods~\cite{Agarwal05, Shashua06}. Guangliang Chen and the
first author~\cite{spectral_theory,spectral_applied} have adapted
the theory developed here and in~\cite{LW-part1,LW-part2} to solve
a problem of multi-way clustering.

\appendix
\section{}

\subsection{Proof of Proposition~\ref{prop:spherical} } \label{app:spherical}

For $\DIm(\HH)>d+1$, the content functions $M_d$ and $M_{d+1}$, and
the norm $\|\cdot\|$ are orthogonally invariant and thus $\pds_0$
and $\gds_0$ are orthogonally invariant. Moreover, in this case,
$M_d$ and $M_{d+1}$ as well as the norm $\|\cdot\|$ scale linearly.
That is, for all $1\leq j\leq d+1$ and $\{\beta_i\}_{i=1}^{d+1}$
such that $\beta_i\neq0,$ where $1 \leq i \leq d+1$:
$$M_d(\beta_1v_1,\ldots,\beta_{j-1}v_{j-1},\beta_{j+1}v_{j+1},\ldots,\beta_{d+1}v_{d+1})=\prod_{i\not=j}|\beta_i|\cdot
M_d\left(v_1,\ldots,v_{j-1},v_{j+1},\ldots,v_{d+1}\right),$$
$$M_{d+1}(\beta_1v_1,\ldots,\beta_{d+1}v_{d+1})=\prod_{i=1}^{d+1}|\beta_i|\cdot
M_{d+1}(v_1,\ldots,v_{d+1})\,,$$ and
$\|\beta_iv_i\|=|\beta_i|\cdot\|v_i\|.$  One can then observe that
both the numerator and denominator of $|\pds_0|$ and $|\gds_0|$
scale similarly and thus the latter functions are invariant under
nonzero dilations. Similarly, the proposition is satisfied when
$\DIm(\HH)=d+1$. \qed

\subsection{On the positivity of the coefficients $\mathbf{\{Q_i\}_{i=1}^{d+1}}$
defined by equation~(\ref{eq:1st-q_i}) } \label{app:sign_q_i} We
show here that the numerators and denominators of the terms $Q_i$,
$1\leq i\leq d+1$, defined by equation~(\ref{eq:1st-q_i}), have the
same signs and thus conclude that these terms are positive.

For $1\leq i\not= j\leq d+1$ we have that
\begin{multline*}
\sign[\gds_{\tilde{u}}(\beta_1 v_1,\ldots,\beta_{j-1}
v_{j-1},0,\beta_{j+1} v_{j+1},\ldots,\beta_{d+1} v_{d+1})]\\=
\sign[\det(\beta_1 v_1-\tilde{u},\ldots,\beta_{j-1}
v_{j-1}-\tilde{u},-\tilde{u},\beta_{j+1}
v_{j+1}-\tilde{u},\ldots,\beta_{d+1} v_{d+1}-\tilde{u})]\\
=-\sign[\det(\beta_1v_1,\ldots,\beta_{j-1}v_{j-1},\tilde{u},\beta_{j+1}
v_{j+1},\ldots,\beta_{d+1} v_{d+1})]\,.\end{multline*}
By the same calculation we also see that
\begin{multline*}\sign[\gds_{\beta_i v_i}(\beta_1
v_1,\ldots,\beta_{j-1} v_{j-1},\tilde{u},\beta_{j+1}
v_{j+1},\ldots,\beta_{i-1} v_{i-1},0,\beta_{i+1}
v_{i+1},\ldots,\beta_{d+1} v_{d+1})] \\= -\sign[\det(\beta_1
v_1,\ldots,\beta_{j-1} v_{j-1},\tilde{u},\beta_{j+1}
v_{j+1},\ldots,\beta_{d+1} v_{d+1})]\,.\end{multline*}
Hence,
\begin{multline*}\sign[\gds_{\tilde{u}}(\beta_1
v_1,\ldots,\beta_{j-1} v_{j-1},0,\beta_{j+1}
v_{j+1},\ldots,\beta_{d+1} v_{d+1})]\\=\sign[\gds_{\beta_i
v_i}(\beta_1 v_1,\ldots,\beta_{j-1} v_{j-1},\tilde{u},\beta_{j+1}
v_{j+1},\ldots,\beta_{i-1} v_{i-1},0,\beta_{j+1}
v_{j+1},\ldots,\beta_{d+1} v_{d+1})]\,,\end{multline*}
and the claim is concluded.

%


\subsection{Proofs of Proposition~\ref{prop:artube} and
Corollary~\ref{cor:arcone}} \label{app:control_msr} We verify here
Proposition~\ref{prop:artube} and Corollary~\ref{cor:arcone}. We
first notice that Corollary~\ref{cor:arcone} is an immediate
consequence of Proposition~\ref{prop:artube} since whenever $x \in
\LL$ we have that
$$
\Cne(\theta,\LL,x) \subseteq \Tbe(\LL,\sin(\theta) \cdot r )\,.
$$

Proposition~\ref{prop:artube} can be concluded from the following
lemma:
\begin{lemma}
\label{lemma:N-m-theta} The set $\Supp \cap \Tbe(\LL,\epsilon \cdot
r )\cap \BB(x,r)$ can be covered by $N$ balls of radius $2\cdot
\sqrt{2} \cdot \epsilon \cdot r$, such that
\begin{equation}\label{eq:N-m-theta}
N \leq \frac{(1+\epsilon)^m}{\epsilon^m} \leq
\frac{2^m}{\epsilon^m}\,.
\end{equation}
\end{lemma}
\begin{proof}
We choose a set $\{y_i\}_{i=1}^N$ in $\Supp \cap \Tbe(\LL,\epsilon
\cdot r )\cap \BB(x,r)$, which is maximally separated by distances
$2\cdot \sqrt{2} \cdot \epsilon \cdot r$. That is,
\begin{equation}\label{eq:yi_contained_in}
\{y_i\}_{i=1}^N \subseteq \Supp \cap \Tbe(\LL,\epsilon \cdot r )\cap
\BB(x,r)\,,
\end{equation}
\begin{equation}\label{eq:sep_yi_cover}
\|y_i-y_j\|> 2\cdot \sqrt{2} \cdot \epsilon \cdot r \,, \ \textup{
for } 1\leq i< j\leq N\,,
\end{equation}
and
\begin{equation}\label{eq:max_sep_yi_cover}
\Supp \cap \Tbe(\LL,\epsilon \cdot r )\cap \BB(x,r) \subseteq
\bigcup_{i=1}^N \BB(y_i, 2 \cdot \sqrt{2} \cdot \epsilon \cdot r)\,.
\end{equation}

We denote $z_i := P_{\LL}(y_i)$, $i=1,\ldots,N$, that is, $z_i$ is
the projection of the point $y_i$ onto the $m$-dimensional affine
plane $\LL$. Equations~(\ref{eq:yi_contained_in})
and~(\ref{eq:sep_yi_cover}) imply that $\{z_i\}_{i=1}^N$ are
separated by distances $2 \cdot \epsilon \cdot r$. Consequently, the
balls $\{\BB(z_i, \epsilon \cdot r)\}_{i=1}^{N}$ are disjoint and
$\{z_i\}_{i=1}^N \subseteq \LL \cap \BB(x,r)$.


We denote by $\mathcal{H}_m$ the $m$-dimensional Hausdorff measure
restricted to $\LL$, and recall that in our case $\mathcal{H}_m$ is
a scaled Lebesgue measure on $\LL$, such that for any ball
$\BB\subseteq\LL$, $\mathcal{H}_m(\BB)=\left(\diam(\BB)\right)^m$.
We thus obtain that
\begin{multline}\label{eq:Hausdorff-covering-number}
N \cdot \left(2 \cdot r\cdot\epsilon\right)^m 
=\sum_{i=1}^{N} \mathcal{H}_m \left(\BB\left(z_i,\epsilon \cdot
r\right)\right)
=\mathcal{H}_m\left(\bigcup_{i=1}^N\BB\left(z_i,\epsilon \cdot r
\right)\right)\\
\leq \mathcal{H}_m \left(\BB\left(x, \left(1+\epsilon\right) \cdot r
\right)\right) =2^m\cdot \left(1+\epsilon\right)^m \cdot r^m\,.
\end{multline}
Equation~(\ref{eq:N-m-theta}) follows directly from
equation~(\ref{eq:Hausdorff-covering-number}) and thus the lemma is
concluded.
\end{proof}
In order to conclude Proposition~\ref{prop:artube} we note that
equation~(\ref{eq:max_sep_yi_cover}) and the definition of an
Ahlfors regular measure imply that
\begin{equation}\label{eq:tube-measure-bound}\mu\left(\Tbe(\LL,\epsilon \cdot
r)\cap\BB(x,r)\right) \ \leq\ \sum_{i=1}^{N}\mu\left(\BB(y_i,2 \cdot
\sqrt{2} \cdot \epsilon \cdot r)\right)\ \leq\ C_{\mu}\cdot N \cdot
2^{\frac{3\gamma}{2}}\cdot\epsilon^{\gamma}\cdot
r^{\gamma}.\end{equation} Then, combining
equations~(\ref{eq:N-m-theta}) and~(\ref{eq:tube-measure-bound}), we
conclude equation~(\ref{eq:artube}).

\section*{Acknowledgement}
We thank Ofer Zeitouni for his helpful comments on an earlier
version of this manuscript. We thank Immo Hahlomaa, Martin
Mohlenkamp and the anonymous reviewers for the very careful
reading of this manuscript and their constructive suggestions. We
thank Paul Nevai for the professional handling of this manuscript.
GL thanks Mark Green and IPAM (UCLA) for inviting him to take part
in their program on multiscale geometry and analysis in high
dimensions. This work has been supported by NSF grant \#0612608.

\bibliographystyle{plain}

\end{document}